\documentclass{amsart}
\usepackage{amsmath,amssymb,enumerate,verbatim}
\usepackage{amsmath, amssymb}
\usepackage{srcltx}
\usepackage{amsthm}
\usepackage{lipsum}
\usepackage[usenames,dvipsnames]{xcolor}
\usepackage{graphicx}
\usepackage{multirow}
\usepackage{booktabs}
\usepackage{mathtools,bm}
\usepackage{blkarray}
\usepackage{tikz} 
\usepackage{multirow}
\usepackage{array,longtable}
\usetikzlibrary{matrix,arrows}
\usepackage{dsfont}
\usepackage{float}
\usepackage{perpage}
\usepackage{mathrsfs}

\usepackage[all]{xy}

\newcommand{\Br}{\overline}

\newcommand{\ep}{{\epsilon}}

\newcommand\la{\lambda}

\newcommand{\bbC}{{\mathbb{C}}}
\newcommand{\bbF}{{\mathbb F}}

\newcommand{\bbR}{{\mathbb R}}

\newcommand{\bbZ}{{\mathbb Z}}

\newcommand{\cal}{\mathcal}

\newcommand{\fk}{\mathfrak}
\newcommand{\ovl}{\overline}
\newcommand{\bb}{\mathbb}

\newcommand{\frakg}{{\mathfrak{g}}}
\newcommand{\frakh}{{\mathfrak{h}}}

\newtheorem{definition}[subsection]{Definition}
\newtheorem{lemma}[subsection]{Lemma}

\newtheorem{prop}[subsection]{Proposition}
\newtheorem{thm}[subsection]{Theorem}
\newtheorem{cor}[subsection]{Corollary}
\newtheorem{theorem}{Theorem}[section]
\theoremstyle{remark}
\newtheorem{remark}[theorem]{Remark}
\newtheorem{notation}[theorem]{Notation}

\newcommand{\calB}{{\mathcal{B}}}

\newcommand{\calF}{{\mathcal{F}}}

\newcommand{\calM}{{\mathcal{M}}}

\newcommand{\calO}{{\mathcal{O}}}

\newcommand{\tu}{\widetilde}

\begin{document}


\title{Lift of the trivial representation to a nonlinear double cover}

\author{Wan-Yu Tsai}

 \thanks{The author is supported by a research grant from the Ministry of Science and Technology of Taiwan.}

\curraddr{
Wan-Yu Tsai \\
Department of Mathematics\\
National Tsing Hua University\\
Hsinchu City \\
Taiwan} 
\email{wytsai@math.nthu.edu.tw}

\begin{abstract}
Let $\tu G$ be the nonlinear double cover of the real points of a connected, simply connected, semisimple complex group. In \cite{Ts}, we introduce a set of genuine small 
representations of $\tu G$ with infinitesimal character $\lambda$, denoted $\prod _\lambda ^s (\tu G)$. In this paper, we show that 
$\prod _{\rho/2} ^s (\tu G)$ is precisely the set of genuine irreducible representations arising from the Kazhdan-Patterson lifting of the trivial representation, when $\tu G$
is simply laced and split. 

\end{abstract}

\maketitle

\section{Introduction}

The  theory  of lifting of characters for nonlinear groups was developed by Kazhdan and Patterson for $GL(n,\bbF)$ in \cite{KP1} and \cite{KP2}. In \cite{AHu}, more specifically, the lifting theory was considered for  $GL(n,\bbR)$: An operator $t_{*}$ which takes representations of $GL(n,\bbR)$ to its nonlinear double cover $\tu{GL}(n,\bbR)$ is defined on the level of global characters, and 
it was shown that if $\pi$ is an irreducible unitary representation of $G$, then $t_* (\pi)$ is either zero or an irreducible unitary representation, up to sign.

Lifting  was studied for more general groups in \cite{AHe}, where  $G$ is assumed to be the real points of a connected, reductive, simply laced complex group, and 
$\tu G$ is an admissible cover of $G$. Identify the kernel of the covering map $p: \tu G\to G$ with $\pm 1$. 
Then the lifting operator was defined by Adams and Herb, stated in the following theorem, which relates genuine characters of $\tu G$ to characters of $G$. By character we mean the 
global character of a representation, viewed as a function on the regular semisimple elements. A genuine representation of $\tu G$ is a representation which does not factor to $G$.
The character of a representation $\pi$ is denoted $\Theta _{\pi}$.

\begin{theorem} \label{t:AHe-main} \emph{(\cite[Theorem 1.6]{AHe})} 
Suppose $G(\bbC)$ is a connected, reductive, simply laced complex group, with real points $G$. We assume that the derived group of $G(\bbC)$ is acceptable (see Definition 2.1 in \cite{AHe}). Suppose 
$\tu G$ is an admissible cover of $G$ (cf. Section 2 in \cite{AHe}). Then we can define the transfer factor 
 $\Delta (h, \widetilde{g})$ satisfying the conditions in Definition \ref{d:Lift}, such that for all stable admissible representation $\pi$ of $G$, 
 \begin{equation} 
 \text{\emph{Lift}}(\Theta_{\pi})(\widetilde{g})=\sum \limits _{\{ h\in G | h^2 = p(\widetilde{g})  \}} \Delta (h, \widetilde{g})\Theta _{\pi} (h)
 \end{equation}
 is the character of a genuine virtual representation $\widetilde{\pi}$ of $\widetilde{G}$, or 0. We say $\widetilde{\pi}$ is the lift of $\pi$ and write $\widetilde{\pi}=$\emph{Lift}$^{\widetilde{G}} _G(\pi)$, where $\Theta _{\widetilde{\pi}}=$\emph{Lift}$_G ^{\widetilde{G}}(\Theta_{\pi})$.

\end{theorem}

Theorem \ref{t:AHe-main} is formally similar to transfer in the setting of endoscopic groups. For example, in \cite{S},  Lift$_G ^{\tu G}$ is analogous to the simplest case of endoscopy: 
 transfer from the quasisplit form $G_{qs}$ of $G(\bbC)$ to $G$. 

\bigskip

In this paper, we adopt Theorem \ref{t:AHe-main}  as the definition of the lifting operator taking representations of $G$ to genuine representations of $\tu G$. 
According to Theorem \ref{t:AHe-main}, Lift$_G ^{\tu G} (\pi)$ is a finite sum of genuine irreducible representations of $\tu G$ with integral coefficients, namely, 
Lift$_G ^{\tu G} (\pi) = \sum _{\tu\pi} a_{\tu \pi}\tu\pi $ with $a_{\tu \pi}\in\bbZ$. We denote Lift$(\pi) = \{\tu\pi \mid a_{\tu \pi} \neq 0\}$ as a finite set. 

In \cite{AHu}, there is a complete discussion of Lift$_G ^{\tu G} (\pi)$
for one-dimensional representations $\pi$ of $GL(n,\bbR)$.
Suppose $G = GL(n, \bbR)$. If $n$ is even, $\tu G$ has a distinguished irreducible representation $T_n$, which comes from the minimal parabolic
 subgroup and contains the pin representation as its lowest $K$-type (see \cite{H} for details). If $n$ is odd, $\tu G$ has two such representations $T_n(\chi_0)$,
 depending on the choice of central character $\chi_0$. Moreover, we write Speh$(k)$ for the $k$th Speh representation of $\tu{GL}(2m,\bbR)$
 with $k=1/2, 1, 3/2, \dots$ (see Section 5 of \cite{AHu}).  Then  Lift($\pi$), with one-dimensional $\pi$ of $G$, is given in the following proposition.
 
 \begin{prop} \emph{(\cite[Proposition 1.5]{AHu})} Suppose $G=GL(2n,\bbR)$.
 \begin{itemize}
 \item[(a)]  If $n$ is even, \emph{Lift}$_G ^{\tu G} (\bbC)=T_n$, and \emph{Lift}$_G ^{\tu G} (sgn)=\text{\emph{Speh}} (1/2)$. 
 \item[(b)]  If $n$ is odd, \emph{Lift}$_G ^{\tu G} (\bbC)=T_n (\chi _0)$, and \emph{Lift}$_G ^{\tu G} (sgn)=0$. 
 \end{itemize}
 
 \end{prop}

Here we consider a similar analysis for other simply laced groups. 
The setting in this paper is as follows. Let $G$ be the real points of a simply connected, semisimple, simply laced complex group, and let $\tu G$ be the nonlinear double cover of $G$.  
Because of the assumptions, the only one-dimensional representation of G is the trivial representation $\bbC$.  
We expect that Lift$_G ^{\tu G}(\bbC)$ gives an interesting class of unitary representations. Recall that in \cite{Ts}, we introduced a set of genuine small representations $\prod _{\rho/2} ^s (\tu G)$, 
which consists of irreducible genuine representations of $\tu G$ with infinitesimal character $\rho/2$ and maximal $\tau$-invariant. Moreover, 
$\prod _{\rho/2} ^s (\tu G)$ is equal to $\prod _{\rho/2} ^{\calO} (\tu G)$, the set of irreducible genuine representations with infinitesimal character $\rho/2$ and attached to a nilpotent orbit $\calO$
(see Table 1 in \cite{Ts} for the list of $\calO$). The following theorem describes the properties that a representation occurring in Lift$(\bbC)$ should possess. 
More precisely,  representations in Lift($\bbC)$ are small representations with infinitesimal character $\rho/2$. This is an analogous result to the case of $\tu{GL}(n,\bbR)$ (cf. \cite[Corollary 4.7]{AHu}).

\begin{thm} \label{t:inclusion}
Suppose $G$ is the real points of a simply connected, semisimple, simply laced complex group, and $\tu G$ is the nonlinear double cover of $G$.  Then 
\begin{equation*}
\textstyle \text{\emph{Lift}}(\bbC) \subseteq  \prod _{\rho/2} ^s (\tu G) =  \prod _{\rho/2} ^{\calO} (\tu G).
\end{equation*}
\end{thm}

We expect that the inclusion in Theorem \ref{t:inclusion} is an equality in certain cases. In this paper, we treat the case of split groups, showing that the inclusion  in 
Theorem \ref{t:inclusion} is an equality
 for split groups. The following is the main theorem of this paper.

\begin{thm} \label{t:main1} \emph{(Main Theorem)}
Suppose $G$ is the real points of a simply connected, semisimple, simply laced complex group, and $\tu G$ is the nonlinear double cover of $G$.  In addition, suppose $G$ is split. 
Then
\begin{equation*}
\textstyle \text{\emph{Lift}}(\bbC) = \prod _{\rho/2} ^s (\tu G) =  \prod _{\rho/2} ^{\calO} (\tu G).
\end{equation*} 
\end{thm}

In \cite{Ts}, we gave an explicit construction of $\prod _{\rho/2} ^s (\tu G)$ when $\tu G$ is split. 
Together with Theorem  \ref{t:main1}, 
this gives a precise description of the irreducible representation that 
occur from lifting of the trivial representation. 

Here is a summary of the contents of the paper. In Section 2, we summarize some notations needed, many of which are from \cite{Ts}. In Section 3, we follow \cite{AHe} to introduce
 the lifting operator and review some of its properties. In Section 4, we focus on lifting the trivial representation. We show that 
representations occurring in Lift$(\bbC)$ must be small with infinitesimal character $\rho/2$. This is Theorem \ref{t:inclusion}. We state the main Theorem at the end of Section 4. 
Sections 5 and 6 are the preparation for proving the main Theorem. We give a brief description of representation theory of nonlinear split groups $\tu{SL}(n,\bbR)$, $\tu{Spin}(n,n)$, and $\tu{E_n}$(split), 
$n=6, 7, 8$ in Section 5. 

The strategy of proving Lift$(\bbC)=\prod _{\rho/2} ^s (\tu G)$ is actually proving Lift$(\bbC)=\prod _{R_D}(\tu G)$. The set $\prod _{R_D}(\tu G)$
is introduced in \cite{Ts}, which gives an explicit description of $\prod _{\rho/2} ^s (\tu G )$ by the theory of Cayley transforms. In \cite{Ts}, it is shown that 
$ \prod _{R_D}(\tu G) \subseteq \prod _{\rho/2} ^s (\tu G)$ and $|\prod _{R_D}(\tu G)|= 1, 4, $ or 16. By a counting argument of Weyl group representations, we show in 
\cite{Ts} that  $\prod _{R_D}(\tu G) = \prod _{\rho/2} ^s (\tu G)$ for type $A$ and $D$. In Section 6 of this paper, we complete the counting argument to show that 
$\prod _{R_D}(\tu G) = \prod _{\rho/2} ^s (\tu G)$ for type E. Then in Section 7, we prove the main Theorem by proving that $\prod _{R_D}(\tu G)\subseteq$Lift$(\bbC)$,
and hence $\prod _{\rho/2} ^s(\tu G)\subseteq$Lift$(\bbC)$, which completes the proof of the main Theorem. 

There are some related questions that we plan to pursue in the future.
First of all, since the lifting operator is defined on the level of global characters, it is natural to ask for the character formula of a small genuine representation which occurs in Lift$(\bbC)$.
When $\tu G= \tu{GL}(n,\bbR)$,  the character formulas of $T_n$ and Speh$(1/2)$ are computed explicitly in \cite[Section 6]{AHu} with the auxiliary of the $t_*$ operator acting on one-dimensional
representations. 

Another interesting question is to study Lift$(\bbC)$ for groups other than split groups. In Corollary \ref{c:nonsplit}, we point out that Lift$(\bbC)=\emptyset$ if $G$ is not in Table 2 of \cite{Ts}.
This means, the non-split simply laced groups, which have nontrivial Lift$(\bbC)$, are the nonlinear double covers of the following groups:
\begin{align*}
A_{2m-1}: \ & \ SU(m,m) \\
A_{2m}: \ & \ SU(m+1,m) \\
D_n:\  & \ Spin(n+1,n-1), \ Spin(n+2,n-2)\\
E_6: \ & E_{6(2)}, \ \text{the quasisplit form of }  E_6(\bbC).
\end{align*}
The groups on this list are all quasisplit, plus $Spin(n+2,n-2)$.
We are hoping to classify  Lift$_G^{\tu G} (\bbC)$ for these groups, to obtain an interesting class of genuine unitary representations of $\tu G$.

\emph{Acknowledgement.} The idea of this paper was originated from the Ph.D. thesis of the author. Part of the contents in this paper were included in her thesis. The author  is grateful to her Ph.D.
advisor Prof. Jeffrey Adams for his constant support for years. 

\section{Some notation}
We always assume that $G$ is a real form of a connected, simply connected, semisimple, simply laced  complex Lie group $G_{\bb C}$, and $\widetilde{G}$ is the unique nonlinear double cover of $G$ 
(see \cite{A2}).  We identify the kernel of the covering map $p: \widetilde{G}\to G$ with $\pm 1$ and write $\widetilde{H}$ for the inverse image in $\widetilde{G}$ of a subgroup $H$ of $G$. 
\medskip

Let $\fk g _\bbR$ be the Lie algebra of $G$ and let 
$\fk g$ be its complexification. Let $K$ denote the maximal compact subgroup of $G$  with corresponding  Cartan involution $\theta$. 
We fix a maximally split $\theta$-stable Cartan subgroup $H^s$ of $G$ with Cartan subalgebra $\fk h_{\bbR} ^s$ (with complexification $\fk h^s$). 
Fix a positive root system $\Delta ^+$ of $\Delta := \Delta (\fk g, \fk h^s)$. Let $W$ be the Weyl group of the root system $\Delta$. 

For a regular element $\lambda$ in $\fk h^*$, where $\fk h^*$ is a Cartan subalgebra of $\fk g$, we will write $\Delta^+(\lambda)$ for the positive root system 
of $\Delta (\fk g, \fk h)$ making $\lambda$ dominant. Let 
$$
R(\lambda) =\{ \alpha \in\Delta (\fk g, \fk h) \mid  \langle \alpha,\lambda \rangle \in \bbZ \}
$$
be the integral root system defined by $\lambda$, and 
$$
R^+(\lambda) =R(\lambda) \cap \Delta ^+ (\lambda), \text{ and } W(\lambda) =W(R(\lambda)),
$$
the positive integral roots and the integral Weyl group. 
\medskip

Let $\cal{HC}(\fk g, K)$ be the set of Harish-Chandra modules and let  $\cal{HC}(\fk g, K)_\lambda \subset \cal{HC}(\fk g, K)$ be the set of Harish-Chandra modules 
with infinitesimal character $\lambda$.  The set of equivalence classes of irreducible admissible representations of $G$, denoted $\widehat G _{adm}$, can be regarded as 
a subset of $\cal{HC}(\fk g, K)$ by sending an irreducible admissible representation to its space of $K$-finite vectors. Similarly, $\widehat G_{adm,\lambda}$ is denoted the set 
of  equivalence classes of irreducible admissible representations of $G$ with infinitesimal character $\lambda$.  The same notions will also be used for the nonlinear group $\tu G$. An
irreducible representation $\pi$ of $\tu G$ is \emph{genuine} if $\pi(-\tu g)=-\pi(\tu g)$ for all $\tu g\in \tu G$. We denote by $\prod_g (\tu G)$ the set of 
equivalence classes of irreducible genuine admissible representations.

\subsection{Regular characters} \label{s:regchar}
  Let $\pi\in \widehat{G}_{adm, \lambda}$, where $\lambda$ is a regular infinitesimal character. Then $\pi$ can be specified by a parameter, which is called a $\lambda$-regular character,  $\gamma=(H,\Gamma,\overline{ \gamma})$, where $H$ is a $\theta$-stable Cartan subgroup of $G$, $\Gamma$ is a character of $H$, and $\overline{\gamma}$ is an element in $\mathfrak{h}^{\ast}$ which defines the same infinitesimal character as $\lambda$, and there are  certain compatibility conditions between $ \overline{\gamma}$  and $\Gamma$ (see Definition 5.3 in \cite{AT}). Write $H=TA$, where $T=H^{\theta}$ and $A$ is the identity component of $\{ h\in H | \theta (h)=h^{-1}\}$. Let $M=Cent_G(A)$. The conditions on $\gamma$ imply that there is a unique relative discrete series representation of $M$, denoted by $\sigma_M$,
with Harish-Chandra parameter $\overline\gamma$, whose lowest $M\cap K$-type has $\Gamma$ as a highest weight. Then we define a parabolic subgroup $P=MN$ such that $\pi=J(\gamma)$, the unique irreducible quotient of a standard representation $I(\gamma)= Ind _P ^G (\sigma _M \otimes 1)$, which is parametrized by $\gamma$ from a $K$-conjugacy class of regular characters for $\lambda$. 

Recall that (see \cite{A1}, for instance) when $\lambda$ is a regular infinitesimal character, $\mathcal{HC}(\mathfrak{g}, K )_{\lambda}$ is parametrized by the set $\mathcal{P}_{\lambda}$ of $K$-conjugacy classes of $\lambda$-regular characters. Furthermore, the following two sets are bases of the Grothendieck group: \begin{center} 
$\{[ J(\gamma)   ] \} _{\gamma\in\mathcal{P}_{\lambda} }$ and $\{[ I(\gamma)   ] \} _{\gamma\in\mathcal{P}_{\lambda} }$. \end{center}

We have the following definition.
\begin{definition} \label{KLVpoly}
Define the change of basis matrix \begin{center} $[J (\delta)]  = \sum \limits _{\gamma \in \mathcal{P}_{\lambda}  }  M(\gamma , \delta)  [I(\gamma)]$  \end{center}
and the inverse matrix \begin{center} $[I   (\delta)]  = \sum \limits _{\gamma \in \mathcal{P}_{\lambda}  }  m(\gamma , \delta)  [J(\gamma)]$.  \end{center}

Here $M(\gamma , \delta)$ and $m(\gamma , \delta)$ are integers and $M(\gamma , \delta) $ are computed by the Kazhdan-Lusztig-Vogan algorithm when $G$ is linear.
\end{definition}

We define the length of a regular character $\gamma$.

\begin{definition} \label{d:length}
Let $\gamma =(H, \Gamma, \ovl\gamma)$ be a regular character of $G$, where $\fk h_{\bbR}$ has Cartan decomposition $ \fk h_{\bbR} = \fk t_{\bbR}+\fk a_{\bbR}$. The 
length of $\gamma$ is  
\begin{equation}
\ell (\gamma)=\frac{1}{2} |\{ \alpha \in \Delta ^+ (\ovl \gamma) \mid \theta (\gamma) \notin \Delta ^+ (\ovl \gamma)\} | +\frac{1}{2}\dim \fk a_{\bbR}. 
\end{equation}

\end{definition}

The definitions  above apply to $\tu G$. 
Let $\tu\pi\in \widehat{\tu G}_{adm,\lambda}$. Then $\tu\pi = J(\gamma)$, the unique irreducible quotient of a standard module $I(\gamma)$, which is parametrized by 
a  \emph{genuine} $\lambda$-regular character  $\gamma=(\tu H, \Gamma, \overline\gamma)$
with $\Gamma$  an irreducible genuine representation of $\widetilde{H} = p^{-1} (H)$. In this case 
$\Gamma$ can be replaced by a character of $Z(\widetilde{H})$, a central character of $\widetilde{H}$, because there is a bijection between 
the set of irreducible genuine representations of $\tu H$ and the set of genuine characters of $Z(\tu H)$ (see Proposition 2.2 \cite{ABPTV}).

\medskip

In this paper, we write $I_G(\gamma)$ and $J_G(\gamma)$ (or $I_{\tu G}(\gamma)$ and $J_{\tu G}(\gamma)$) to emphasize the representations of $G$ (or $\tu G$, 
respectively) when it is needed.
\medskip

For nonlinear covering groups, the "extended integral length" $\tu \ell ^I$ is defined for a genuine regular character $\gamma$ in \cite[Definition 6.7]{RT}.
This is the correct notion needed for the induction in the Kazhdan-Lusztig algorithm for nonlinear groups. Note the $\ell (\gamma)=\tu \ell ^I (\gamma)$ if $\tu G$ is simply laced, which 
is the case we will be focusing on throughout the paper. 
\medskip

If $G$ is assumed to be simply laced,  we have the following facts.
\begin{prop} \label{p:center} 
Assume that $G$ is a real form of a connected, simply connected, semisimple, simply laced  complex Lie group $G_{\bb C}$, and $\widetilde{G}$ is the unique nonlinear double cover of $G$. In addition, 
suppose $G$ is simply laced. Let $H$ be a Cartan subgroup of $G$, and let $H^0$ be the identity component of $H$. Then\\
(1) \emph{(\cite{AHe}, Proposition 4.7)} $Z(\widetilde{H})=Z(\widetilde{G})\widetilde{H^0}$. In particular, a genuine character of $Z(\widetilde{H})$ is determined by its restriction to $Z(\widetilde{G})$ and its differential;\\
(2) \emph{(\cite{AT}, Proposition 5.5)} A genuine regular character $\gamma= (\widetilde{H}, \Gamma,\overline{ \gamma})$ of $\widetilde{G}$ is determined by $\overline{\gamma}$ and the restriction of $\Gamma$ to $Z(\widetilde{G})$, and so is $\widetilde{\pi}=J(\gamma)$.
\end{prop}

We also recall the notions of the cross action and Cayley transform defined on  the regular characters (see \cite{AT}). 
\begin{itemize}
\item \textbf{Cross action.}  The cross action of $W$ on regular characters is denoted $w\times \gamma$, where $w\in W$, $\gamma$ is a $\lambda$-regular character. If $w\in W(\lambda)$, then 
$w\times \gamma$ is again a $\lambda$-regular character.
\item\textbf{Cayley transforms and inverse Cayley transforms.}   
Let $\gamma$ be a $\lambda$-regular character. Suppose $\alpha$ is a noncompact imaginary root for $\gamma$. Then $c^{\alpha}(\gamma)=\{ \gamma ^{\alpha}\}$ 
or $\{ \gamma_+^\alpha, \gamma_-^\alpha \}$ denotes the Cayley transform of $\gamma$ by $\alpha$, where $ \gamma^\alpha, \gamma_{\pm}^\alpha$ are again $\lambda$-regular characters.
Suppose $\alpha$ is a real root for $\gamma$. Then $c_{\alpha}(\gamma)=\{ \gamma _{\alpha}\}$ 
or $\{ \gamma^+_\alpha, \gamma^-_\alpha \}$ denotes the inverse Cayley transform of $\gamma$ by $\alpha$, where $ \gamma_\alpha, \gamma^{\pm}_\alpha$ are again $\lambda$-regular characters.

\end{itemize}

 \bigskip

 We recall some more notation from \cite{Ts}. 
 \begin{itemize}
\item  The set of  irreducible genuine \textbf{small} representations of $ \tu G$   with infinitesimal  character  $\lambda$  is 
$$
\textstyle\prod_{\lambda} ^s (\tu G)  = \{ \pi \in \widehat G _{\text{adm}, \lambda} \mid \pi \text{ has maximal $\tau$-invariant} \},
$$ 
where 
$\lambda =\rho/2$ if $G$ is simply laced or of type $G_2$, and see Table 1 in \cite{Ts} for $\lambda$ in other cases. 

Moreover,  $\pi \in \prod_{\lambda} ^s (\tu G) $ if and only if 
 $\pi$ is attached to the nilpotent orbit $\calO$ (see Table 1 in \cite{Ts}).

\medskip

Furthermore, suppose that $G$ is simply laced and split.
\medskip

\item $\mathcal{SH}_{\rho/2}$ denotes the set of irreducible quotients of pseudospherical principal series representations (see \cite{ABPTV} for details) with infinitesimal character $\rho/2$;  
$\mathcal{SH}_{\rho/2 ,\chi}$ denotes the subset of $\cal{SH}_{\rho/2}$ containing representations with central character $\chi$. 

 \item $\prod_{R_D}(\tu G)$ denotes the set of representations constructed by applying Cayley transforms  
 to $Sh\in\cal{SH}_{\rho/2}$ through subsets of simple roots from $R_D$ (see Table 3 in \cite{Ts}); 
 $\prod_{R_D}(\tu G)_{\chi}$  denotes the subset of $\prod_{R_D}(\tu G)$ containing representations with central character $\chi$.
 
 
 \end{itemize}

\section{Kazhdan-Patterson Lifting }
The materials in this section can be found in \cite{AHe}. The assumptions for $G, G_\bbC, \tu G$ are as in Section 2. 
 Recall that we identify the kernel of the covering map $p: \tu G\to G$ with $\{\pm 1\}$.
Let $G'$ be the set of regular semisimple elements of $G$.

\begin{definition} \label{d:stable}
\emph{Let  $\pi \in \widehat{G}_{adm}$, with character $\Theta _{\pi}$.  We say $\pi$ and $\Theta _{\pi}$ are \textit{stable} if $\Theta _{\pi}$ is invariant under conjugation of $G_{\mathbb{C}}$, that is, $\Theta _{\pi} (g)=\Theta _{\pi}(g')$ if $g, g'\in G'$ and $g'=xgx^{-1}$ for some $x\in G_\bbC$.}

\end{definition}

Suppose $H$ is a  Cartan subgroup of $G$ and $\Phi ^+$ is a set of positive roots of $H$ in $G$. For $h\in H$ we have the Weyl denominator \begin{center} $|D(h)|^{\frac{1}{2}}=|\prod \limits _{\alpha \in \Phi ^+} (1-\alpha ^{-1}(h))| |e^{\rho }(h) | $.
\end{center}  

\begin{definition} \label{d:Lift}  \emph{ (see \cite {AHe}) Suppose $\pi\in \widehat{G}_{adm}$ and $\pi$  is stable. For  $\widetilde{g}\in \widetilde{G'}$, define \begin{center} Lift$_G ^{\widetilde{G}}(\Theta_{\pi})(\widetilde{g})=\sum \limits _{\{ h\in G | h^2 = p(\widetilde{g})  \}} \Delta (h, \widetilde{g})\Theta _{\pi} (h)$.  \end{center} Here $\Delta (h,\widetilde{g})$ is a certain function on $G'\times \widetilde{G'}$ satisfying the following conditions: 
\begin{center}
$\Delta (h,\widetilde{g}) =0$ unless $h^2=p(\widetilde{g})$\\
$|\Delta (h, \widetilde{g}) | = |D(h)|^{\frac{1}{2}} /|D (\widetilde{g}) |^{\frac{1}{2}}$\\
$\Delta (xhx^{-1}, \widetilde{x}\widetilde{g}\widetilde{x}^{-1}) =\Delta (h, \widetilde{g}) \quad (\widetilde{x}\in \widetilde{G}, x=p(\widetilde{x}))$\\
$\Delta (h, -\widetilde{g})=-\Delta (h, \widetilde{g})$
\end{center}
By Section 5 in \cite{AHe}, since $G_{\mathbb{C}}$ is simply connected and semisimple, the function $ \Delta$ is canonical.}
\end{definition}
    
     The following theorem is a special case of the main theorem of \cite{AHe}. Since $G_{\mathbb{C}}$ is simply connected and semisimple, a simplified version of Section 5 in \cite{AHe} applies.
     
\begin{theorem}  \label{Liftvirtual}   
Assume the setting in the beginning of this section. Then there is a canonical function (see Section 5 in \cite{AHe}) $\Delta (h, \widetilde{g})$ satisfying the conditions in Definition \ref{d:Lift}, such that for all stable admissible representation $\pi$ of $G$, \begin{center}  \emph{Lift}$(\Theta_{\pi})(\widetilde{g})=\sum \limits _{\{ h\in G | h^2 = p(\widetilde{g})  \}} \Delta (h, \widetilde{g})\Theta _{\pi} (h)$ \end{center} is the character of a genuine virtual representation $\widetilde{\pi}$ of $\widetilde{G}$, or 0. We say $\widetilde{\pi}$ is the lift of $\pi$ and write $\widetilde{\pi}=$\emph{Lift}$^{\widetilde{G}} _G(\pi)$, where $\Theta _{\widetilde{\pi}}=$\emph{Lift}$_G ^{\widetilde{G}}(\Theta_{\pi})$.
\end{theorem}    

Because of this theorem, for a stable admissible representation $\pi$ of $G$, we will denote Lift$ _{G} ^{\widetilde{G}}(\pi)$ as a set as follows. If   Lift$ _{G} ^{\widetilde{G}}(\pi)=\sum \limits _{\widetilde{\pi}} a_{\widetilde{\pi}} \widetilde{\pi}$, for $a_{\widetilde{\pi}}\in \mathbb{Z}$ and $\widetilde{\pi}\in \prod _g(\widetilde{G})$, the set of  irreducible genuine admissible representations of    $\widetilde{G}$, then we define 
 \begin{equation} 
\text{Lift} (\pi)=\{  \widetilde{ \pi} \in \textstyle\prod _g (\widetilde{G})\thinspace | \thinspace a_{\widetilde{\pi}} \neq 0 \}
\end{equation} 
as a set, and this is a finite set of irreducible genuine representations due to Theorem \ref{Liftvirtual}. 

\subsection{Lifts of regular characters and standard modules}

Lifting of regular characters is defined in \cite{AHe}. We summarize some important facts here.

\begin{definition}\label{d:modchar}
\emph{(cf. \cite[Definition 16.2]{AHe}) A }modified regular character
\emph{is a triple $\gamma =(H, \Gamma,\ovl \gamma)$, where $H$ is a $\theta$-stable Cartan subgroup of $G$, $\Gamma$ is a character of $H$ and $\ovl \gamma\in \fk h^*$.
Write $H=TA$ and $M=Cent_G(A)$ with complexified Lie algebra $\fk m$. Let $$\Phi _i=\Delta (\fk m, \fk h),$$  and assume 
$$\langle \lambda, \alpha^\vee \rangle \in \bbR^\times \text{ for all } \alpha\in \Phi _i.$$ 
Let $\Phi_i^+(\ovl\gamma)=\{\alpha \in \Phi_i \mid \langle \ovl\gamma,\alpha^\vee \rangle >0\}$, $\rho_i (\gamma) =\frac{1}{2}\sum_{\alpha\in \Phi_i^+(\ovl\gamma) } \alpha.$  
We assume $d\Gamma =\lambda -\rho_i(\ovl\gamma)$.
}    
\end{definition}

The set of modified regular characters for $G$ is denoted $CD(G)$. If $H$ is fixed, let $CD(G,H)$ be the set of modified regular characters $(H,\Gamma,\ovl\gamma)$.
Associated to $\gamma\in CD(G)$ is a relative discrete series representation $\pi_M (\gamma)$, with Harish-Chandra
parameter $\ovl\gamma$ and central character $\Gamma |_{Z(M)}$. The fact is that $CD(G)$ is in one-to-one bijection with the set of regular characters of $G$ introduced in 
\ref{s:regchar}.  Thus, we will also write $I(\gamma)=I_G(\gamma):=Ind_P ^G (\pi_M(\gamma))$, where $P=MN$ is any parabolic subgroup containing $M$, to be the standard representation of $G$ 
parametrized by $\gamma\in CD(G)$.   Let $\Theta _G(\gamma)$ denote the character of $I_G(\gamma)$.

Let $H$ be a Cartan subgroup of $G$.
Note $W_i=W(\Phi_i)$ acts on $H$ and this induces an action on $CD(G,H)$: $w(H,\Gamma,\ovl\gamma)=(H, w\Gamma,w\ovl\gamma)$ where $w\Gamma(h)=\Gamma(w^{-1}h)$.
Note that $I(w\gamma)$ and $I(\gamma)$ have the same infinitesimal and central characters.

\begin{lemma} \emph{(cf. \cite[Lemma 16.12]{AHe})}
Suppose $\gamma=(H,\Gamma,\ovl\gamma) \in CD(G,H)$. We define 
\begin{equation}
I_G ^{st} (\gamma) =\sum_{w\in W(M,H) \backslash W_i } I_G(w\gamma)
\end{equation} 
with character $\Theta^{st}_G(\gamma)$.
Then $\Theta_G^{st} (\gamma)$ is a stable character (see Definition \ref{d:stable}).
\end{lemma}

The definition of modified regular character extends naturally to $\tu G$ just like that of regular character. 

A  \emph{genuine modified regular character} of $\tu G$ is $\tu \gamma=(\tu H, \tu \Gamma, \ovl\gamma)$ where $\tu H=p^{-1}(H)$ is a Cartan subgroup of $\tu G$, 
$\tu \Gamma$ is a character of $Z(\tu H)$ and $\ovl\gamma\in \fk h^*.$ The conditions for $\tu \gamma$ to satisfy are the same as those for  a modified regular character 
(see Definition \ref{d:modchar}). 

Let $CD_g(\tu G)$ be the set of genuine modified regular characters of $\tu G$ and let $CD_g(\tu G, \tu H)$ be the subset with given Cartan subgroup $\tu H$. 
 Write $I_{\tu G} (\gamma)$ to be the standard representation of $\tu G$ 
parametrized by $\tu \gamma\in CD_g(\tu G)$. 

\medskip

We define lifting on the level of  modified regular characters. 

Given a Cartan subgroup $H$.  Let $\phi$ be the map on $G$ defined by $\phi(g) = g^2$ and write $\phi_H$ for the restriction of $\phi$ on $H$. For $\tu g\in \tu H$, let 
\begin{equation}
X(H,\tu g)= \{h\in H\mid \phi (h)= p(\tu g)\} 
\end{equation}
and 
\begin{equation}
\tu S =\{ \tu g \in \tu H\mid p(\tu g) \in \phi(H)\} = \{ \tu g\in \tu H\mid X( H, \tu g)\neq \emptyset\}.
\end{equation}

Fix $\chi\in \widehat H$ and $\tu\chi\in \prod_g(Z(\tu H))$ to be the lifting data (see \cite[Definition 7.1]{AHe}). For a character $\psi$ of $H$, we have (cf. Section 10 of \cite{AHe})
\begin{equation}
\text{Lift}_H ^{\tu H} (\psi) (\tu g) = \begin{cases}
0 & \text{ if } \tu g\notin \tu S,\\
c^{-1} |\text{Ker} (\phi_H)| \tu{\psi}_0 (\tu g) & \text{ if } \tu g \in \tu S,
\end{cases}
\end{equation}
where $c\in\bbZ$ is given in \cite[Definition 8.1]{AHe}, and 
$\tu \psi _0$ is defined to be
\begin{equation} \label{e:psi0}
\tu\psi _0 (\tu g) =\tu\chi(\tu g) \psi (h) \chi(h)^{-1} 
\end{equation}
with some choice of $h\in X(H,\tu g)$. 
Note that this is independent of the choice of $h$ and hence $\tu \psi _0$ is a well-defined genuine character of $\tu S$.

Let 
\begin{equation} \label{e:Xpsi0}
\textstyle \prod_g (Z(\tu H), \tu \psi _0) =\{\tu\psi \in \prod_g (Z(\tu H)) \mid \tu \psi |_{\tu S} =\tu \psi _0 \}.
\end{equation}

By Proposition 2.2 in  \cite{ABPTV}, there is a bijection between $\prod_g (Z(\tu H))$ and $\prod _g (\tu H)$. We denote this bijection by $\tu \psi \to \tu \tau (\tu \psi)$.
The following proposition summarizes how Lift$_H ^{\tu H} (\psi)$ is defined.

\begin{prop}\label{p:liftchar}
 \emph{(\cite[Proposition 10.11]{AHe})} 
Fix $\tu\chi \in \prod _g (Z(\tu H))$ and $\chi \in \widehat H$, and use them to define \emph{Lift}$_H ^{\tu H}$. 
Fix $\psi \in \widehat H$. Then \emph{Lift}$_H ^{\tu H} (\psi) \neq 0$ if and only if $\chi (t) =\psi (t)$ for all $t\in Ker (\phi _H)$. 
Assume this holds. Define $\tu\psi_0$ as in (\ref{e:psi0})  and $\prod _g (Z(\tu H), \tu\psi_0)$ as in (\ref{e:Xpsi0}). Then 
\begin{equation}
\text{\emph{Lift}}_H ^{\tu H}  (\psi) =\sum _{\tu \psi \in \prod_g (Z(\tu G) ,\tu\psi_0 )} \tu \tau (\tu \psi). 
\end{equation}
This is an identity of genuine representations of $\tu H$; the right hand side is the direct sum of $|p (Z(\tu H) ) / \phi_H(H)|$ irreducible representations. 	
The differentials satisfy 
\begin{equation}
d\tu \psi =\frac{1}{2} (d\psi - d\mu)
\end{equation}
where $d\mu =d\chi -2d\tu \chi\in\fk h^*$.
\end{prop}

Now we are ready to define Lift$_G ^{\tu G}$ of modified regular characters.

\begin{definition} \label{d:liftchar} \emph{(\cite[Definition 17.5]{AHe})
Fix  $\chi, \tu\chi$ as in Proposition \ref{p:liftchar}. Let $\gamma=(H,\Gamma,\ovl\gamma)\in CD(G)$.  If $\Gamma | _{Ker(\phi_H)}\neq \chi | _{Ker (\phi_H)}$ then Lift$_H ^{\tu H} (\Gamma)=0$,
and we define Lift$_G ^{\tu G} (\gamma)=0$ and $\Theta_{\tu G} (\text{Lift}_G ^{\tu G} (\gamma))=0$.
Otherwise write 
\begin{equation}
\text{Lift}_H ^{\tu H} (\Gamma) = \sum ^n _{i=1} \tu \tau(\tu \Gamma_i)
\end{equation}
where $n=|p(Z(\tu H))/\phi (H) |$, and each $\tu \Gamma_i$ is a genuine one-dimensional representation of $Z(\tu H)$ (see Proposition \ref{p:liftchar}).
For $1\le i\le n$ let $\tu \gamma_i =(\tu H, \tu\Gamma_i, \frac{1}{2} (\ovl\gamma- d\mu))$ (with $\mu$ defined in Proposition \ref{p:liftchar}); this is a genuine regular character of $\tu G$.
Define 
\begin{equation}
\text{Lift}_G ^{\tu G} (\gamma) =\{ \tu \gamma_1, \dots, \tu\gamma_n\}
\end{equation}
and 
\begin{equation}
\Theta_{\tu G} (\text{Lift}_G ^{\tu G} (\gamma) ) =\sum _{i=1} ^n \Theta_{\tu G}(\tu \gamma _i). 
\end{equation}
}
\end{definition}

Furthermore, lifting of a stable sum of  standard modules can also be defined. 
Let $\gamma=(H,\Gamma,\ovl\gamma)\in CD(G,H)$. Recall that $W_i=W(\Phi_i)$ acts on $CD(G,H)$.

\begin{theorem} \label{Liftstable} \emph{(\cite[Corollary 19.8]{AHe})}
Suppose $\gamma=(H,\Gamma,\ovl\gamma)\in CD(G,H)$.
Let $\{ \widetilde{\gamma_1}, \cdots, \widetilde{\gamma_n} \}$ be the set of constituents of \emph{Lift}$^{\widetilde{G}} _G (w\gamma)$ as $w$ runs over $W_i$, considered without multiplicity. Then
 \begin{equation}
 \text{\emph{Lift}}^{\widetilde{G}} _G (I_G ^{st}(\gamma))=C(H) \sum\limits _{i=1} ^n I_{\widetilde{G}} (\widetilde{\gamma _i}), 
 \end{equation} 
 where $C(H) =c(H)/c(H_s)$, $c(H)=|H_2 ^0| |H/Z_0 (H)| ^{\frac{1}{2}}$, $H_s$ is the maximally split Cartan subgroup of $G$, $H_2 ^0$ is the subgroup of elements of order 2 in the identity component of $H$, $Z_0(H)=p(Z(\widetilde{H}))$. Note that all the constituents have distinct central characters, so are a fortiori distinct,  and
that $C(H)$ is normalized so that $C(H_s)=1$.

\end{theorem}

\section{Lift of the trivial representation}
In this section,
the assumptions for $G$ and $\tu G$ are the same as in Section 2. 
Denote by $\bbC$ the trivial representation of $G$. In this section, more specifically,  we will study Lift$(\bbC)$.
\medskip

When $G=GL(n, \mathbb{R})$,  and $\pi$ is a one-dimensional representation of $G$, it is shown in Section 4 of \cite{AHu} that 
any $\tu \pi \in$ Lift$(\pi)$ has infinitesimal character $\rho/2$ and maximal $\tau$-invariant, assuming that $\tu\pi\neq 0$. 
We will generalize this result to various groups $G$.  Note that the proof is parallel to that of the 
case of $GL(n,\bbR)$, with a little modification.

\begin{lemma} \label{l:Liftinfchar}
Retain the setting for $G$ in the beginning of this section. 
Let $\pi$ be a stable admissible virtual representation of $G$ with infinitesimal character $\lambda$. Assume that \emph{Lift}$(\pi) \neq \emptyset$, then every $\widetilde{\pi} \in$\emph{Lift}$(\pi)$ has infinitesimal character $\lambda /2$.
\begin{proof}
We just need to show the case for standard modules since standard modules span virtual modules. 

  Let $I ^{st}_G(\gamma)$ be a stable sum of standard modules, parametrized by a modified regular character $\gamma=(H, \Gamma, \lambda)\in CD(G)$. 
   By Definition \ref{d:liftchar}, if Lift$_H ^{\widetilde{H}}(\Gamma)\neq 0$, then  Lift$^{\widetilde{G}} _G (\gamma)=\{   \widetilde{
  \gamma   _1}  , \cdots , \widetilde{  \gamma   _n}   \}$, where each $\widetilde{
  \gamma   _i}  =(\widetilde{H}, \widetilde{\Gamma_i} , \displaystyle\frac{1}{2} (\lambda -\mu) )$ is a genuine regular character of $\widetilde{G}$. It turns out that each $\widetilde{\gamma _i}$ has infinitesimal character $\lambda /2$ since $G_{\mathbb{C}}$ is simple and simply connected, lifting is canonical and $\mu=0$ (see Section 5 in \cite{AHe}). Now by Theorem \ref{Liftstable}, Lift$^{\widetilde{G} }_G ( I ^{st}_G(\gamma))$ is a sum of $I_{\widetilde{G}} (\widetilde{\gamma_i})$ and hence the infinitesimal character of Lift$^{\widetilde{G} }_G ( I ^{st}_G(\gamma))$  is $\lambda/2$. 
\end{proof}
\end{lemma}

Next, we  show that the Lift operator commutes with coherent continuation.

\begin{definition}\emph{ (cf. Definition 7.2.5 in \cite{Vgr}) \label{d:cohfam}
Let $H\subseteq G$ be a Cartan subgroup. The weight lattice in $\widehat H$ is the subgroup $\Lambda$ of $\widehat H$ consisting of weights of finite dimensional representations 
of $G$.  Fix $\lambda\in \fk h^*$. A \textbf{coherent family} based at $\lambda$ is a collection of virtual modules $\{\pi(\lambda +\mu) \mid \mu\in \Lambda\}$ such that for all $\mu
\in\Lambda$,
\begin{itemize}
\item[(a)] $\pi(\lambda +\mu)$ has infinitesimal character $\lambda +d\mu$, and 
\item[(b)] for any finite dimensional representation $F$ of $G$, 
$$
\pi(\lambda +\mu) \otimes F \simeq \sum\limits _{\xi\in\triangle (F)} \pi(\lambda +(\mu+\xi)),
$$
where $\triangle (F)\subset \Lambda$ is the set of weights of $F$. 
\end{itemize} }
\end{definition}

Since every finite dimensional representation of $\tu G$ factors to $G$ when the rank of $G$ is bigger than 1, Definition \ref{d:cohfam} applies when $G$ is replaced with $\tu G$. 

\medskip

Fix $\lambda\in\fk h^*$ and a coherent family $\{ \pi(\lambda + \mu)\mid \mu\in\Lambda \}$ for $G$. 
For $\mu\in\Lambda,$ define 
$$
\tu \pi (\lambda/2 +\mu ) = \text{Lift}_{G}^{\tu G} ( \pi(\lambda + 2\mu). 
$$
 
\begin{prop} \label{p:cohfam2}
The collection of virtual representations $\{   \tu \pi (\lambda/2 +\mu ) \mid \mu\in\Lambda   \}$ is a coherent family.
\end{prop}

We need the following Lemma before giving the proof. The next Lemma for $GL(n,\bbR)$ can be found in \cite{AHu},

\begin{lemma} \label{l:Liftfin}
Retain the setting of $G$ in the beginning of the section. Let $F$ be a finite dimensional representation of G. Then there exists a 
virtual finite dimensional representation $A \ast F$ satisfying 
\begin{equation} \label{e:Fchar}
\Theta_{A\ast F}(g) = \Theta_F(g)^2
\end{equation}
for all $g\in G$. Equivalently, for all Cartan subgroups $H$ of $G$,
\begin{equation} \label{e:wt}
\triangle (A\ast F) =2\triangle (F), 
\end{equation}
where $\triangle (F ) , \triangle (A\ast F)\subset \widehat H$ denote the weights of $F$ and $A\ast F$ respectively.
 \end{lemma}

See Section II.7 in \cite{BT} for the case of compact groups G. The result can be extended to general Lie groups.

\begin{proof}[proof of Proposition \ref{p:cohfam2}]
By Lemma \ref{l:Liftinfchar}, $\tu \pi (\lambda /2 +\mu)$ is a virtual representation with infinitesimal character $\lambda/2 +d\mu$. We need to show that for any finite dimensional 
representation $F$ of $\tu G$ and $\mu\in\Lambda$, we have 
\begin{equation}\label{e:cohfam3}
\sum\limits _{\xi\in\triangle(F)} \tu \pi(\lambda/2 + \mu +\xi) = \tu \pi (\lambda/2 +\mu) \otimes F.
\end{equation}
Note that $F$ is a nongenuine finite dimensional representation of $\tu G$, that is, for $\tu g \in \tu G$,
\begin{equation}
\Theta_{F}(\tu g) = \Theta_F(p(\tu g)). 
\end{equation}
We now compute the character of both sides of (\ref{e:cohfam3}) at a regular semisimple element $\tu g\in\tu G$, which means that 
$\tu g$ is a preimage of a regular semisimple element of $g\in G$. 
\medskip

The following sums over $h$ are over $\{h\in G\mid h^2 = p(\tu g)\}$.

\begin{align*}
\sum\limits_{\xi\in\triangle(F)} \Theta_{\tu\pi (\lambda/2 +\mu +\xi)} (\tu g )
&= \sum\limits_{\xi\in\triangle(F)} \sum\limits_h   \Delta(h,\tu g)  \Theta_{\pi (\lambda +2\mu +2\xi)} (h) & (\text{Definition } \ref{d:Lift})\\
&= \sum\limits_h   \sum\limits_{\xi\in\triangle( A\ast F)}  \Delta(h,\tu g)  \Theta_{\pi (\lambda +2\mu +\xi)} (h) & ((\ref{e:wt})) \\ 
&= \sum\limits_h  \Delta(h,\tu g)   \Theta_{\pi (\lambda +2\mu )} (h) \Theta _{A\ast F} (h)  & (\text{Definition \ref{d:cohfam}(b)} )\\
&= \sum\limits_h  \Delta(h,\tu g)   \Theta_{\pi (\lambda +2\mu )} (h) \Theta _{ F} (p(\tu g))& ((\ref{e:Fchar}))\\ 
&= \Theta_{\tu\pi (\lambda/2 +\mu)} (\tu g ) \Theta_F (\tu g) & (\text{Definition } \ref{d:Lift})
\end{align*}
This completes the proof.
\end{proof}

\begin{theorem} \label{t:Liftsmall}
Let $\widetilde{\pi}\in$\emph{Lift}$(\bbC)$. Then $\widetilde{\pi}$ has infinitesimal character $\rho/2$ and maximal $\tau$-invariant.
\begin{proof}
The first assertion of the theorem is a corollary of Lemma \ref{l:Liftinfchar} since the infinitesimal character of $\bbC$ is $\rho$. 

To show that $\widetilde{\pi}$ has maximal $\tau$-invariant, we will show that $\psi_{\alpha}(\widetilde{\pi})=0$ for all $\alpha\in \prod(\rho/2)$, where $\psi_{\alpha}$ is the Zuckerman translation functor pushing to the $\alpha$-wall.  Suppose $\lambda$ is singular and $\rho/2-\lambda$ is a sum of roots. Since the trivial representation $\bbC$ is finite-dimensional, 
the translation of $\bbC$ to infinitesimal character is zero; by Proposition \ref{p:cohfam2}, the same holds for the translation of $\tu \pi$ to $\lambda$.
\end{proof}
\end{theorem}

\begin{cor} \label{c:liftsubset}
Retain the setting in the beginning of the section for $G$ and $\widetilde{G}$. Then
\begin{equation}
\textstyle \text{\emph{Lift}}(\bbC )\subseteq \prod _{\rho/2} ^s (\widetilde{G})
= \prod _{\rho/2} ^{\mathcal{O}} (\widetilde{G}).
\end{equation}
\end{cor}

This corollary limits the irreducible representations $\tu \pi$ that can contribute to the lift of the trivial representation: $\tu \pi$ must have infinitesimal character
$\rho/2$ and maximal $\tau$-invariant.

\begin{cor} \label{c:nonsplit}
\emph{Lift}$(\bbC) =\emptyset$ if $G$ is not in Table 2 in \cite{Ts}. Therefore, if \emph{Lift}$_G ^{\widetilde{G}} (\bbC) \neq 0$ then $G$ is quasisplit or $G = Spin (n+2,n-2).$
\end{cor}

We expect that the inclusion in Corollary \ref{c:liftsubset} is an equality for many groups $\tu G$. It turns out that this is true when $\tu G$ is split.  
The following is the main theorem of the paper. We will complete the proof of it in Section \ref{s:main}.

\begin{theorem} \emph{(\textbf{Main Theorem})}
Suppose that  $G$ is the split real form of a simply laced, connected, simply connected complex group, with $\tu G$ to be the nonlinear double cover of $G$.Then 
$$
\textstyle
 \emph{\text{Lift}}(\bbC )= \prod _{\rho/2} ^s (\widetilde{G}),
$$
that is,  the set of genuine irreducible representations occurring in the lifting of  $\bbC$ is precisely the set of genuine small representations of $\tu G$ with 
infinitesimal character $\rho/2$.

\end{theorem}

\section{Representations  of Split groups $\tu G$}

We will focus on simply laced split groups from now on. Therefore, $\tu G =\tu{SL}(n,\bbR)$ for type $A_{n-1}$, $\tu G=\tu{Spin}(n,n)$ for type $D_n$, 
or $\tu G$ is the nonlinear double cover for split real form of type $E_6, E_7,$ or $E_8$ for type $E$. For later use, we list all   
$\rho/2$-regular parameters of $\tu G$.  

\subsection{Type $A_{n-1}$} \label{s:typeA}
$\tu G =\tu{SL}(n,\bbR)$. 
The Cartan subgroups of $\tu G$ can be enumerated as $\tu H_i$, where $i$ is the real rank of $\tu H _i$. More 
precisely,  when $n=2p$, we have $p-1\le i\le n-1$, and 
\begin{align*}
H_{p-1}&\cong  S^1 \times  (\bbC^\times)^{p-1}  \\
H_i &\cong   (\bbC ^\times)^{n-1-i} \times (\bbR ^{\times}) ^{2i-n+1}  \quad \text{ for }  p \le i\le n-1;
\end{align*}
when $n=2p+1$, we have $p\le i \le n-1$, and 
\begin{align*}
H_{p}&\cong  S^1 \times  (\bbC^\times)^{p-1} \times \bbR^\times  \\
H_i &\cong   (\bbC ^\times)^{n-1-i} \times (\bbR ^{\times}) ^{2i-n+1}  \quad \text{ for }  p +1 \le i\le n-1.
\end{align*}

\begin{lemma} \emph{(cf. \cite{ABPTV}) } \label{l:Z-A}
\begin{itemize}
\item[(a)] For $n=2p$, $Z(\tu G) =\bbZ_2\times \bbZ_2$.
\item[(b)] For $n=2p+1$, $Z(\tu G) =\bbZ _2$.
\end{itemize}
\end{lemma}

Fix $\Delta ^+=\{ e_i-e_j \mid i<j\}$ for all parameters. Denote $\Delta ^+_1 = \Delta _{\rho/2} ^+ $ to be the positive integral roots for $\rho/2$, and 
$\Delta _{1/2}$ to be the half-integral roots for $\rho/2$. Therefore, we have $\Delta ^+ =\Delta ^+ _1 \cup \Delta ^+ _{1/2}$.
More precisely, $$\Delta ^+ _1= \{ e_i-e_j \mid i< j, \ i-j\in 2\bbZ\}, \ 
\Delta ^+ _{1/2}= \{ e_i-e_j \mid i< j, \ i-j\in 2\bbZ +1\}.$$

Due to Lemma \ref{l:Z-A}, $\cal{SH}_{\rho/2}$ consists of  2 genuine representations  with different central characters when $n=2p$ and it consists of  1 representation when $n=2p+1$.
Write $\gamma_{Sh}$ for the parameter(s) of $Sh\in\cal{SH}_{\rho/2 }$. All parameters of $\tu G$ can be obtained from $\gamma_{Sh}$ by (inverse) Cayley transforms through real
i roots  (for $Sh$) in $\Delta^+ _{1/2}$. We write parameters specified by Cartan subgroups in this way:
\begin{align*}
\tu H_{n-1}: \  &  \gamma_{Sh}\\
\tu H_{n-1-j}: \  & c_{\{ \alpha_1, \dots, \alpha_j\} } (\gamma_{Sh}), \text{ where }  \{\alpha_1,\dots,\alpha_j\} \subset \Delta^+ _{1/2}\ \text{ is an orthogonal set}. 
\end{align*}

Fix a central character $\chi$ of $\tu G$, and take $Sh\in\cal{SH}_{\rho/2}$ with central character $\chi$.  A $\rho/2$-regular character $\gamma = c_{\{e_{i_1} -e_{j_1}, \cdots, e_{i_m} -e_{j_m} \}} (\gamma_{Sh})$
 with central character $\chi$
can be parametrized by a set of pairs of numbers 
$$
\{ \{i_1,j_1\}, \{ i_2,j_2\}, \dots, \{i_m, j_m\} \}
$$
with $i_k\in \{1, 3, \dots, 2p-1\}$ if $n=2p$, (or $i_k\in \{1, 3, \dots, 2p-1, 2p+1\}$ if $n=2p+1$) and 
$j_k\in \{2, 4, \dots, 2p\}$. 

\medskip

We write  $\gamma = c_{\{e_{i_1} -e_{j_1}, \cdots, e_{i_m} -e_{j_m} \}} (\gamma_{Sh}) = \gamma ( \{i_1,j_1\}, \{ i_2,j_2\}, \dots, \{i_m, j_m\}) $.
We write $\gamma ( S_1)\subseteq \gamma (S_2)$ if $S_1\subseteq S_2$.

The corresponding standard module and irreducible quotients are denoted 
$$I( \{i_1,j_1\}, \{ i_2,j_2\}, \dots, \{i_m, j_m\}  )$$
and
$$J(\{i_1,j_1\}, \{ i_2,j_2\}, \dots, \{i_m, j_m\}  ),$$
respectively. 

In particular, 
\begin{equation}
\textstyle \prod _{R_D} (\tu G)  _\chi=\begin{cases}
 \{J(\emptyset) \} & \text{ if } n=2p+1,\\
 \{J( \emptyset) , J(\{1,2\}, \{3,4\}, \dots, \{2p-1,2p\}) \} & \text{ if } n=2p.
\end{cases}
\end{equation}

Recalling the length of $\gamma$ defined in \ref{d:length} and letting $\ell_s =\ell (\gamma (\emptyset))$,  we have that 
$$
\ell (\gamma (\{ i_1, j_1\}, \dots,\{ i_m, j_m\}    )) = \ell _s -\sum_{k=1} ^m |i_k-j_k|,
$$
if we can arrange the $i_k$ and $i_k$ so that $i_1<j_1<\cdots<i_m<j_m$.
Note that $\ell (\gamma(S_1)) > \ell (\gamma (S_2))$ if $S_1\subset S_2$.

\subsection{Type $D_n$} \label{ss:D}

We refer the reader to \cite[Section 6]{BTs} for more details about the representation theory of $\tu{Spin}(n,n)$. 

As in the case of type $A$, a $\rho/2$-regular character $\gamma$ with central character $\chi$ can be obtained from some $Sh\in\cal{SH}_{\rho/2}$ by 
a series of Cayley transforms, i.e. 
$\gamma=c_{\alpha_1,\dots,\alpha_m}(Sh)$, where $\{\alpha_1,\dots,\alpha_m\}$ is a set of strongly orthogonal roots.

Such $\gamma$ can be parametrized by a set of the form
\begin{equation} \label{e:D-param}
\{ \{ i_1, i_2,\dots, i_r, j_1, j_2,\dots, j_r\} ,  \{\ep_{r+1}i_{r+1}  ,   \ep_{r+1}j_{r+1}  \},  \{\ep_{r+2}i_{r+2}  ,   \ep_{r+2}j_{r+2}  \} , \dots, \{ \ep_m i_m, \ep_m j_m\} \}
\end{equation}
with $\ep_k=\pm 1$, $i_k\in \{1, 3, \dots, 2p-1\}$ if $n=2p$, (or $i_k\in \{1, 3, \dots, 2p-1, 2p+1\}$ if $n=2p+1$) and 
$j_k\in \{2, 4, \dots, 2p\}$.  If $\gamma$ is a regular character parametrized by (\ref{e:D-param}), then it is denoted
$$
\gamma(\{ \{ i_1, i_2,\dots, i_r, j_1, j_2,\dots, j_r\} ,  \{\ep_{r+1}i_{r+1}  ,   \ep_{r+1}j_{r+1}  \},  \{\ep_{r+2}i_{r+2}  ,   \ep_{r+2}j_{r+2}  \} , \dots, \{ \ep_m i_m, \ep_m j_m\} \} );
$$
it means $e_{i_k} \pm e_{i_l}$, $e_{j_k} \pm e_{j_l}$ are compact imaginary for $\gamma$, 
$e_{i_k} \pm e_{j_l}$ are noncompact imaginary for $\gamma$ for $1\le k\neq l \le r$, and $e_{i_k} -\ep_k e_{jk}$ is noncompact imaginary for $\gamma$ for $r+1\le k \le m$. 

\medskip
We write $\gamma ( S_1)\subseteq \gamma (S_2)$ if $S_1\subseteq S_2$. Note that $\ell (\gamma(S_1)) > \ell (\gamma (S_2))$ if $S_1\subset S_2$.

The corresponding standard module and irreducible quotients are denoted 
$$I(\{ \{ i_1, i_2,\dots, i_r, j_1, j_2,\dots, j_r\} ,  \{\ep_{r+1}i_{r+1}  ,   \ep_{r+1}j_{r+1}  \},  \{\ep_{r+2}i_{r+2}  ,   \ep_{r+2}j_{r+2}  \} , \dots, \{ \ep_m i_m, \ep_m j_m\} \} )$$
and
$$J(\{ \{ i_1, i_2,\dots, i_r, j_1, j_2,\dots, j_r\} ,  \{\ep_{r+1}i_{r+1}  ,   \ep_{r+1}j_{r+1}  \},  \{\ep_{r+2}i_{r+2}  ,   \ep_{r+2}j_{r+2}  \} , \dots, \{ \ep_m i_m, \ep_m j_m\} \} ),$$
respectively. 

In particular, 
\begin{equation}
\textstyle \prod _{R_D} (\tu G)_\chi =\begin{cases}
 \{J(\emptyset) ,  J( \{ 2p-1, 2p \} , \{ -(2p-1) , -2p\} ) \} & \text{ if } n=2p+1,\\
 \{J( \emptyset) , J( \{ 2p-1, 2p \} , \{ -(2p-1) , -2p\} ),   & \text{ if } n=2p.\\
 J(\{1,2\}, \{3,4\}, \dots \{2p-3, 2p-2\}, \ep\{2p-1,2p\}) \}, \ep=1 \text{ or }  -1 & 
\end{cases}
\end{equation}


\subsection{Type $E$}
Let $G_\bbC$ be a complex algebraic group of type $E_n$, $n=6, 7, 8$, with Lie algebra $\fk g$. 
Let $G$ be the split real form of $G_\bbC$ with Lie algebra $\fk g _\bbR$. Let $\theta$ denote the Cartan involution (corresponding to $G$).

\subsection{Cartan subalgebra in $\fk g_ \bbR$} \label{s:E-CSA}
Conjugacy classes of Cartan subalgebras have the following representatives:
\begin{equation}
\begin{aligned}
E_6: & \ \fk h^{r, m , s} _\bbR =\{ (x_1,x_2,x_3,x_4,x_5, -x_6,-x_6,x_6)\}\\
E_7: &\  \fk h^{r, m , s} _\bbR=\{ (x_1,x_2,x_3,x_4,x_5, x_6,-x_7,x_7)\}\\
E_8: &\  \fk h^{r, m , s} _\bbR =\{ (x_1,x_2,x_3,x_4,x_5, x_6,x_7,x_8)\}
\end{aligned}
\end{equation}
with $\theta$ given in the  tables below.

\begin{align*}
E_6 :\  &  \fk h _{\bbR}  ^{r,m,s} \ \        &&\hskip-2mm   \theta   \ \ &&\hskip-2mm   c_{r,m,s}\\     
           &  \fk h _{\bbR} ^{2,2,0} \ \            &&\hskip-2mm   (-Id)\circ s_{12} s_{\ovl{12}} s_{34} s_{\ovl{34}}   \ \ &&\hskip-2mm  c_{12} c_{\ovl{12}} c_{34} c_{\ovl{34}} \\                  
              & \fk h_{\bbR} ^{0,3,0} \ \     &&\hskip-2mm   (-Id)\circ s_{12} s_{\ovl{12}} s_{34}   \ \ &&\hskip-2mm  c_{12} c_{\ovl{12}} c_{34}  \\ 
           & \fk h_{\bbR} ^{0,2,2} \ \     &&\hskip-2mm   (-Id)\circ s_{12} s_{\ovl{12}}    \ \ &&\hskip-2mm  c_{12} c_{\ovl{12}}   \\ 
             & \fk h_{\bbR} ^{0,1,4} \ \     &&\hskip-2mm   (-Id)\circ s_{12}    \ \ &&\hskip-2mm  c_{12}   \\ 
             & \fk h_{\bbR} ^{0,0,6} \ \     &&\hskip-2mm   -Id  \ \ &&\hskip-2mm  Id  \\
\end{align*}
\begin{align*}
  E_7 :\  &  \fk h _{\bbR} ^{r,m,s} \ \        &&\hskip-2mm   \theta   \ \ &&\hskip-2mm   c_{r,m,s}\\   
  &  \fk h_{\bbR}  ^{7,0,0}\ \   &&\hskip-2mm    (-Id)\circ s_{1 2}s_{  \overline{12}} s_{34}s_{\Br{34}} s_{56} s_{\Br{56}} s_{78}  \ \ &&\hskip-2mm   c_{1 2}c_{  \overline{12}} c_{34}c_{\Br{34}} c_{56} c_{\Br{56}} c_{78}\\     
            &  \fk h_{\bbR}  ^{5,1,0}\ \   &&\hskip-2mm    (-Id)\circ s_{1 2}s_{  \overline{12}} s_{34}s_{\Br{34}} s_{56}  s_{78}  \ \ &&\hskip-2mm   c_{1 2}c_{  \overline{12}} c_{34}c_{\Br{34}} c_{56}   c_{78} \\     
             &  \fk h_{\bbR}  ^{3,2,0}\ \   &&\hskip-2mm    (-Id)\circ s_{1 2}s_{  \overline{12}} s_{34}s_{\Br{34}}   s_{78}  \ \ &&\hskip-2mm  c_{1 2}c_{  \overline{12}} c_{34}c_{\Br{34}}   c_{78} \\   
              &  \fk h_{\bbR}  ^{1,3,0}\ \   &&\hskip-2mm    (-Id)\circ s_{1 2}s_{  \overline{12}} s_{34}   s_{78}  \ \ &&\hskip-2mm c_{1 2}c_{  \overline{12}} c_{34}   c_{78}   \\  
               &  \fk h_{\bbR}  ^{2,2,1}\ \   &&\hskip-2mm    (-Id)\circ s_{1 2}s_{  \overline{12}} s_{34}s_{\Br{34}}   \ \ &&\hskip-2mm   c_{1 2}c_{  \overline{12}} c_{34}c_{\Br{34}} \\      
                 &  \fk h_{\bbR}  ^{1,2,2}\ \   &&\hskip-2mm    (-Id)\circ s_{1 2}s_{  \overline{12}} s_{78}  \ \ &&\hskip-2mm  c_{1 2}c_{  \overline{12}} c_{78} \\  
                   &  \fk h_{\bbR}  ^{0,3,1}\ \   &&\hskip-2mm    (-Id)\circ s_{1 2}s_{  \overline{12}} s_{34}   \ \ &&\hskip-2mm c_{1 2}c_{  \overline{12}} s_{34} \\ 
                      &  \fk h _{\bbR} ^{0,2,3}\ \   &&\hskip-2mm    (-Id)\circ s_{1 2}s_{  \overline{12}}    \ \ &&\hskip-2mm   c_{1 2}c_{  \overline{12}}\\  
                      &  \fk h _{\bbR} ^{0,1,5}\ \   &&\hskip-2mm    (-Id)\circ s_{1 2}   \ \ &&\hskip-2mm c_{1 2}   \\  
                     &  \fk h _{\bbR} ^{0,0,7}\ \   &&\hskip-2mm    -Id   \ \ &&\hskip-2mm  Id 
  \end{align*}
  \begin{align*}
  E_8 :\  &   \fk h_{\bbR}  ^{r,m,s} \ \        &&\hskip-2mm   \theta   \ \ &&\hskip-2mm   c_{r,m,s}\\   
&   \fk h _{\bbR} ^{8,0,0}\ \   &&\hskip-2mm    (-Id)\circ s_{1 2}s_{  \overline{12}} s_{34}s_{\Br{34}} s_{56} s_{\Br{56}} s_{78} s_{\Br{78}}  \ \ &&\hskip-2mm  
  c_{1 2}c_{  \overline{12}} c_{34}c_{\Br{34}} c_{56} c_{\Br{56}} c_{78} c_{\Br{78}} \\     
            &  \fk h _{\bbR} ^{6,1,0}\ \   &&\hskip-2mm    (-Id)\circ s_{1 2}s_{  \overline{12}} s_{34}s_{\Br{34}} s_{56} s_{\Br{56}}  s_{78}  \ \ &&\hskip-2mm  
            c_{1 2}c_{  \overline{12}} c_{34}c_{\Br{34}} c_{56} c_{\Br{56}}  c_{78} \\     
             &  \fk h_{\bbR}  ^{4,2,0}\ \   &&\hskip-2mm    (-Id)\circ s_{1 2}s_{  \overline{12}} s_{34}s_{\Br{34}}  s_{56}  s_{78}  \ \ &&\hskip-2mm c_{1 2}c_{  \overline{12}} c_{34}c_{\Br{34}}  c_{56}  c_{78}  \\   
              &  \fk h _{\bbR} ^{2,3,0}\ \   &&\hskip-2mm    (-Id)\circ  s_{1 2}s_{  \overline{12}} s_{34}  s_{\Br{34}} s_{56}  \ \ &&\hskip-2mm  c_{1 2}c_{  \overline{12}} c_{34}  c_{\Br{34}} c_{56} \\  
               &  \fk h _{\bbR} ^{0,4,0}\ \   &&\hskip-2mm    (-Id)\circ s_{1 2}s_{  \overline{12}} s_{34} s_{56}  \ \ &&\hskip-2mm   c_{1 2}c_{  \overline{12}} c_{34} c_{56}  \\      
                 &  \fk h_{\bbR}  ^{2,2,2}\ \   &&\hskip-2mm    (-Id)\circ s_{1 2}s_{  \overline{12}} s_{34} s_{\Br{34}}  \ \ &&\hskip-2mm  c_{1 2}c_{  \overline{12}} c_{34} c_{\Br{34}}  \\  
                   &  \fk h _{\bbR} ^{0,3,2}\ \   &&\hskip-2mm    (-Id)\circ s_{1 2}s_{  \overline{12}} s_{34}   \ \ &&\hskip-2mm  c_{1 2}c_{  \overline{12}} c_{34}   \\ 
                      &  \fk h_{\bbR}  ^{0,2,4}\ \   &&\hskip-2mm    (-Id)\circ s_{1 2}s_{  \overline{12}}    \ \ &&\hskip-2mm  c_{1 2}c_{  \overline{12}}  \\  
                      &  \fk h _{\bbR} ^{0,1,6}\ \   &&\hskip-2mm    (-Id)\circ s_{1 2}   \ \ &&\hskip-2mm  c_{1 2} \\  
                       &  \fk h_{\bbR}  ^{0,0,7}\ \   &&\hskip-2mm    -Id   \ \ &&\hskip-2mm  Id 
\end{align*}

\textit{Remark and Notation.}
We write $\fk h^{r,m,s}$ for the complexification of $\fk h^{r,m,s} _\bbR$. Note that $\fk h_{\bbR} ^{0,0,n}$ is the split Cartan subalgebra. In the second column of the table,
$s_{ij}( x_i) =x_j$, $s _{\ovl{ij}} = -x_j$. In the last column,
 $c_{ij}$ denotes the (inverse) Cayley transform through the root $e_i-e_j$, and $c_{\Br{ij}}$ denotes the (inverse) Cayley transform through the root $e_i+e_j$. 
The notation $c_{m,r,s}$ relates $\fk h ^{r,m,s}$ to $\fk h ^{0,0,n}$ as follows:
$$
\fk h_{\bbR}^{r,m,s} = c_{m,r,s} (\fk h _\bbR ^{0,0,n}),
$$ 
where $c_{m,s,r}$ is a product of (inverse) Cayley transforms $c_\alpha$ shown in the tables. 
The Haase diagrams of Cartan subalgebras are exhibited in Figure 1. The numbers on the left of the figure indicate the real rank of each Cartan subalgebra. 
\begin{figure}[h]
\caption{}
\begin{center}
     \begin{picture}(500,300)(0,-7)
                  \put(88,265){$E_6$}
                   \put(198,265){$E_7$}
                  \put(313,265){$E_8$}
	\put(90,180){\circle{4}}
	\put(96,180){$\fk h _\bbR ^{2,2,0}$ }
	\put(90, 175){\line(0,-1){20}}
	\put(90, 150 ){\circle{4}}
	\put(96, 150 ){$\fk h _\bbR ^{0,3,0}$}
	\put(90,145 ){\line(0, -1){20}}
	\put(90,120){\circle{4}}
	\put(96,120){ $\fk h _\bbR ^{0,2,2}$}
          \put(90,115){\line(0,-1){20}}	
          \put(90,90){\circle{4}}
           \put(96,90){ $\fk h _\bbR ^{0,1,4}$}
          \put(90, 85  ){\line(0,-1){20}}
           \put(90,60){\circle{4}}
             \put(96,60){$\fk h _\bbR ^{0,0,6}$}
             \put(15,237){0}
             \put(15,207){1}
          \put(15,177){2}
          \put(15, 147){3}
           \put(15,117){4}
           \put(15, 87){5}
           \put(15,57){6}
          \put(15,27){7}
              \put(15,-3){8}
              
              \put(200,240){\circle{4}}
	   \put(206,240){ $\fk h _\bbR ^{7,0,0}$}
	\put(200, 235){\line(0,-1){20}}
	\put(200, 210 ){\circle{4}}
	\put(206, 210 ){$\fk h _\bbR ^{5,1,0}$ }
	\put(200,205 ){\line(0, -1){20}}
	\put(200,180){\circle{4}}
	\put(206,180){ $\fk h _\bbR ^{3,2,0}$}
       \put(198,175){\line(-1, -1){23}}
      \put(202 ,175 ){\line(1, -1){23}}
	
	    \put(174,147){\line(2, -1){50}}
	    
    \put(172,115){\line(1, -1){23}}
        \put(228 ,115 ){\line(-1, -1){23.5}}

   	\put(230,150){\circle{4}}
	   	\put(234,150){ $\fk h _\bbR ^{2,2,1}$ }

		\put(230,145 ){\line(0, -1){20}}
	\put(170,150){\circle{4}}
	\put(156,153){$\fk h _\bbR ^{1,3,0}$   }
		\put(170,145 ){\line(0, -1){20}}
	
   	\put(230,120){\circle{4}}
	   	\put(234,124){ $\fk h _\bbR ^{0,3,1}$ }
		
	\put(170,120){\circle{4}}
	\put(148,123){$\fk h _\bbR ^{1,2,2}$  }
          \put(200,90){\circle{4}}
          \put(197,100){$\fk h _\bbR ^{0,2,3}$  }
    
          \put(200, 85  ){\line(0,-1){20}}
           \put(200,60){\circle{4}}
                    \put(206,60){ $\fk h _\bbR ^{0,1,5}$     }
          \put(200, 55  ){\line(0,-1){20}}
           \put(200,30){\circle{4}}
                    \put(206,30){ $\fk h _\bbR ^{0,0,7}$ }

  \put(313,145){\line(-1, -1){23.5}}
        \put(317 ,145 ){\line(1, -1){23.5}}

    \put(287,115){\line(1, -1){23.5}}
        \put(343 ,115 ){\line(-1, -1){23.5}}

    \put(315,240){\circle{4}}
        \put(321,240){$\fk h _\bbR ^{8,0,0}$  }
    
	\put(315, 235){\line(0,-1){20}}
	\put(315, 210 ){\circle{4}}
		\put(321, 210 ){$\fk h _\bbR ^{6,1,0}$ }
	\put(315,205 ){\line(0, -1){20}}
	\put(315,180){\circle{4}}
	\put(320,180){ $\fk h _\bbR ^{4,2,0}$  }
	\put(315,175 ){\line(0, -1){20}}
	\put(315,150){\circle{4}}
		\put(321,150){$\fk h _\bbR ^{2,3,0}$  }
	
   	\put(345,120){\circle{4}}
	   	\put(350,120){  $\fk h _\bbR ^{2,2,2}$  } 
	\put(285,120){\circle{4}}
	\put(266,120){ $\fk h _\bbR ^{0,4,0}$  }
	
          \put(315,90){\circle{4}}
                \put(315,80){  $\fk h _\bbR ^{0,3,2}$  }
          
          \put(315, 84 ){\line(0,-1){20}}
           \put(315,60){\circle{4}}
                  \put(320,60){   $\fk h _\bbR ^{0,2,4}$}
           
          \put(315, 55  ){\line(0,-1){20}}
           \put(315,30){\circle{4}}
          \put(320,30){  $\fk h _\bbR ^{0,1,6}$ }
        \put(315, 25  ){\line(0,-1){20}}
           \put(315,0){\circle{4}}
           \put(320,0){ $\fk h _\bbR ^{0,0,8}$}             
                   \end{picture} 
\end{center}
\end{figure}

\begin{lemma}
The $\fk h _\bbR ^{r,m,s}$ listed above exhaust the $G$-conjugacy classes of Cartan subalgebras in $\fk g_\bbR$. 
\end{lemma}

\subsection{Center of $\tu G$}

Let $A=MA^0$ be the split torus of $G$ with identity component $A^0$. Then $\tu A =\tu M A^0$ is a two-fold cover of $A$. By Lemma 3.12
in \cite{ABPTV}, $Z(\tu G)=Z(\tu M)$. Furthermore, $p(Z(\tu M)) \cong [2P^{\vee} \cap R^\vee]/2R^\vee$, where $P^\vee$ and $R^\vee$ are the coweight lattice and coroot lattice,
respectively. Consequently, we have the following lemma.

\begin{lemma} \label{l:Z-E} \emph{(cf. \cite{ABPTV})} 
\begin{itemize}
\item[(1)] In the case of type $E_6$, $Z(\tu G)\cong \bbZ_2$.
\item[(2)] In the case of type $E_7$, $Z(\tu G)\cong \bbZ_2\times \bbZ_2$.
\item[(3)] In the case of type $E_8$, $Z(\tu G)\cong \bbZ_2$.
\end{itemize}
\begin{proof}
By the comments above this lemma, it suffices to compute $Z(\tu M)$.
In the cases of $E_6$ and $E_8$, $|[2P^{\vee} \cap R^\vee]/2R^\vee| =1$.
In the case of $E_7$, the representatives for $[2P^{\vee} \cap R^\vee]/2R^\vee$ can be picked to be
$\{0, (0,\dots, 0, -1,1)\}$ and hence $| [2P^{\vee} \cap R^\vee]/2R^\vee | =2$. This gives the desired results of this lemma.
\end{proof}
\end{lemma}

\subsection{Cartan subgroups of $\tu G$ }
Recall that in \ref{s:E-CSA}, $\fk h _{\bbR} ^{r,m,s}$ denotes a Cartan subalgebra of $G=E_n$(split). Let $H^{r,m,s}$ denote the centralizer of $\fk h_\bbR ^{r,m,s}$ in $G$ of $\fk h _\bbR ^{r,m,s}$.
We have 
$$
H^{r,m,s}\cong (S^1)^r\times (\bbC ^{\times}) ^m \times (\bbR ^\times) ^s. 
$$
The centralizer of $\fk h_\bbR ^{r,m,s}$ in $\tu G$ gives a Caratn subgroup of $\tu G$. It turns out that such a Cartan subgroup is $\tu H^{r,m,s} = p^{-1} (H^{r,m,s})$. 

\subsection{Regular characters}\label{s:char-E}
 We use coordinates from \cite{Bou} for the root system in type $E$. 
\begin{notation} \label{n:E}
A root of the form $\frac{1}{2} (\pm 1, \dots,\pm 1)$ will be denoted 
\begin{center} $\beta _{i_1\cdots i_k}$, \ $1\le i_1 < \cdots < i_k \le 8$, \end{center}
 if $-1$ occurs in the $i_j$-th coordinate. For example, $\beta _{234567} = \frac{1}{2} (1,-1,\dots,-1,1)$ is one of the simple roots in type $E$.  The simple reflections through 
 the root $e_i-e_j$ and $e_i+e_j$ will be denoted $s_{ij}$ and $s_{\ovl{ij}}$, respectively; the simple reflection through the root $\beta _{i_1,\dots,i_k}$ will be denoted $r_{i_1 i_2\cdots i_k}$.
 
\end{notation}

Fix $\Pi$ to be the set of simple roots of type $E$, given as follows:
\begin{align*}
E_6 : \  & \  \Pi = \{ \alpha_1 = \beta_{234567},\  \alpha_2 = e_1+e_2,  \  \alpha_3= e_2-e_1, \  \alpha_4 = e_3-e_2 , \\
 & \ \alpha_5= e_4-e_3, \ \alpha_6 = e_5-e_4\},  \\  
E_7: \  &  \  \Pi = \{ \alpha_1 = \beta_{234567},\  \alpha_2 = e_1+e_2,  \  \alpha_3= e_2-e_1, \  \alpha_4 = e_3-e_2 , \\
 & \ \alpha_5= e_4-e_3, \ \alpha_6 = e_5-e_4, \  \alpha_7 = e_6-e_5\},  \\  
E_8: \  &  \  \Pi = \{ \alpha_1 = \beta_{234567},\  \alpha_2 = e_1+e_2,  \  \alpha_3= e_2-e_1, \  \alpha_4 = e_3-e_2 , \\
 & \ \alpha_5= e_4-e_3, \ \alpha_6 = e_5-e_4, \  \alpha_7 = e_6-e_5, \ \alpha_8= e_7-e_6 \}.
\end{align*}

We summarize 
$\rho/2$, $\Delta (\rho/2),$ and $\Pi(\rho/2)$ for each type as follows. Again, we use coordinates from \cite{Bou}. 

\begin{align*}
           & \rho/2    \   && \Delta (\rho/2) \  &&  \Pi (\rho/2) \\     
E_6 :\  &  \frac{1}{2}(0,1,2,3,4,-4,-4,4)  \    &&   A_1\times A_5 \    &&    e_2+e_4,  -e_1+e_3, -e_3+e_5,  e_1+e_3, \\
              & &&&& \beta_{134567}, -e_2+e_4    \\ 
  E_7 :\  &  \frac{1}{2} (0,1,2,3,4,5,-\frac{17}{2}, \frac{17}{2})      \ &&  A_7 \  &&  -e_1+e_3, -e_3+e_5,  e_1+e_3,  \beta_{134567},  \\ 
       & &&&&  -e_2+e_4, -e_4+e_6, e_2+e_4 \\
   E_8 :\  &   \frac{1}{2} (0,1,2,3,4,5,6,23)  \ &&   D_8  \ &&   e_2+e_4,  -e_4+e_6, -e_2+e_4, \beta_{134567}, \\
       & &&&&  e_1+e_3, -e_3+e_5, -e_5+e_7, -e_1+e_3
  \end{align*}

Fix a set of positive roots $\Delta ^+$ corresponding to $\Pi$ for all parameters. Let $\Delta _1 ^+=\Delta (\rho/2)$, which is formed of $\Pi (\rho/2)$ and 
$\Delta ^+ _{1/2} =\{ \alpha\in \Delta ^+\mid \langle  \rho/2, \alpha^\vee\rangle \in \bbZ+1/2\}$.  

Due to Lemma \ref{l:Z-E}, $\cal{SH}_{\rho/2}$ consists of  one genuine representation   when $\tu G = \tu{E _6}$(split) and  $\tu{E _8}$(split) and it consists of  2 genuine  representations when 
$\tu G = \tu{E _7}$(split).
Write $\gamma_{Sh}$ for the parameter(s) of $Sh\in\cal{SH}_{\rho/2 }$. All parameters of $\tu G$ can be obtained from $\gamma_{Sh}$ by (inverse) Cayley transforms through real
roots  (for $Sh$) in $\Delta ^+ _{1/2}$. We write parameters specified by Cartan subgroups in this way:
\begin{align*}
\tu H^{0,0,n}: \  &  \gamma_{Sh}\\
\tu H^{r,m,s}: \  & c_{\{ \alpha_1, \dots, \alpha_j\} } (\gamma_{Sh}), \text{ where }  \{\alpha_1,\dots,\alpha_j\} \subset \Delta ^+ _{1/2}\ \text{ is an orthogonal set}. 
\end{align*}

Moreover, if $\alpha_1,\dots,\alpha_j$ are all of the form $e_i\pm e_k$, then the parameter $c_{\{ \alpha_1, \dots, \alpha_j\} } (\gamma_{Sh})$ can be written as the form 
$j_k\in \{2, 4, \dots, 2p\}$.  If $\gamma$ is a regular character parametrized by (\ref{e:D-param}), then it is denoted
$$
\gamma(\{ \{ i_1, i_2,\dots, i_r, j_1, j_2,\dots, j_r\} ,  \{\ep_{r+1}i_{r+1}  ,   \ep_{r+1}j_{r+1}  \},  \{\ep_{r+2}i_{r+2}  ,   \ep_{r+2}j_{r+2}  \} , \dots, \{ \ep_m i_m, \ep_m j_m\} \} );
$$
as in \ref{ss:D} for type $D$; similar notation also applies to the corresponding standard modules and irreducible quotients. 

In particular, fixing a central character $\chi$ of $\tu G$,

\begin{equation}
\textstyle \prod _{R_D} (\tu G)_\chi =\begin{cases}
 \{J(\emptyset) \} & \text{ if }  \tu G = \tu{E_6}(\text{split}),\\
 \{J( \emptyset) , J( \{ -1, -2\}, \{3, 4\}, \{5, 6 \} )\}  & \text{ if }  \tu G = \tu{E_7}(\text{split}),\\
 \{J(\emptyset) \} & \text{ if }  \tu G = \tu{E_8}(\text{split}),\\
\end{cases}
\end{equation}

\section{Counting Small representations for type $E$}

We focus on the group $G$ which is a semisimple, simply laced, connected, split real group and assume $\tu G$ is the nonlinear double cover of $G$. We will count the number of
small representations in $\prod_{\rho/2} ^s(\tu G)$.  The number of representations attached to $\calO$ with infinitesimal character $\rho/2$ equals the multiplicity of the sign representation of $W(\rho/2)$ in the coherent continuation representation.  For type $A_{n-1}$ and $D_n$, this has been done in \cite{Ts}. We will treat the case of type $E$ in this section.

We first recall the coherent continuation action for the group $\tu G$.

\subsection{Coherent continuation representations for $\tu G$}
We take a look at the cross action again. As in \cite{RT},  given the infinitesimal character $\lambda$, we define a family of 
infinitesimal characters $\calF (\la)$ including $\la$. Every $\la '\in \calF (\lambda)$ can be indexed by some $w\in W/W(\lambda)$. 
Write $\calB _{\lambda ', \chi}$ for the set of equivalence classes of standard representation parameters with infinitesimal character $\lambda ' \in\calF(\lambda)$ and a fixed central character $\chi$ of $\tu{G}$, and $$\calB := \coprod \limits _{\lambda '\in \calF(\lambda) , \chi \in   \prod_g (Z(\tu{G}) )    } \calB_{\lambda ' , \chi}.$$ 

As we will see later the coherent continuation action is closely related to the cross action. 
We may use $\calF(\lambda)$ to define the cross action of $W$ on $\calB$, denoted $w\times \gamma$ for $w\in W$, $\gamma\in\calB$, as shown in \cite{RT}. In fact, fixing an infinitesimal character $\lambda ' \in \calF(\lambda)$ and a central character $\chi$, $W(\lambda)$ acts on $\calB_{\lambda ,\chi}$ by the cross action.

\smallskip

On the other hand, there is a coherent continuation action of $W(\lambda)$ on $ \bbZ[\calB_{\lambda,\chi}]$. The formulas of the coherent continuation action can be derived from those of the Hecke operators as we shall see.  We set $\calM = \bbZ [u^{ \frac{1}{2}}, u^{-\frac{1}{2}} ] [\calB]$.
We fix the abstract infinitesimal character $\la_a\in \calF(\la)$ corresponding to the   positive root system $\Delta ^+:=\Delta ^+ _a (\frakg, \frakh^a)$ (where $\frakh^a$ is an abstract Cartan subalgebra of $\frakg$) and the set of simple roots $\Pi_a \subset \Delta ^+ _a$. For $s=s_{\alpha}$ with $\alpha\in \Pi _a$, the action of $T_s$ on $\gamma\in \calM$ is defined in Section 9 of \cite{RT}. 
\smallskip

We  also consider $\Delta (\la)$, the integral root system for $\la$ and the integral Weyl group $W(\la)$. 
Choose $\Pi(\lambda)$ to be a set of simple roots for $\Delta(\lambda)$. 
\medskip

Given $\alpha\in \Pi(\lambda)$, we will need to decompose $s_\alpha$ into a product of simple reflections. This can 
be easily done. In fact, we need the decomposition in the following Lemma.

\begin{lemma} \label{l:chamber}
Suppose $\Delta$ is a simply laced root system with $\Pi\subset \Delta ^+$ a choice of simple roots. Let $\alpha\in\Delta$.
Suppose $s_\alpha$ is decomposed into a product of simple reflections  
$$s_\alpha = s_{\alpha_{2l +1} } \cdots s_{\alpha_2} s_{\alpha_1}$$ 
with   $\alpha_j\in \Pi$ and $\alpha_j = \alpha_{2l+2-j}$ for $1\le j \le l$.   Define $\beta_k = s_{\beta_{k-1}}\cdots s_{\beta_1}$ for 
$1\le k \le 2l +1$.
Then 
\begin{equation}\label{e:simple}
s_{\alpha}= s_{\beta_{2l +1} }\cdots s_{\beta_2}s_{\beta_1}
\end{equation}
with each $\beta_k$ simple for the chamber of $s_{\beta_{k-1} }\cdots s_{\beta_1} (\rho)$.
\begin{proof}
This can be proved by induction on $l$. 
\end{proof}
\end{lemma}

Given a root $\beta,$ we write $T_\beta$ for the Hecke operator $T_{s_{\beta}}$ for simplicity. Using the notation in Lemma \ref{l:chamber}, 
\begin{equation} \label{hecke}
\begin{aligned}
T_{\alpha} (\gamma) &=T_{\beta_{2l +1}} \cdots  T_{\beta_2}  T_{\beta_1} (\gamma) \\
&= p_{2l +1}(u) \cdots p_1 (u)  s_{\beta_{2l +1}}\times (s _{ \beta _{2l }}  \times \cdots (s_{\beta_1 }\times \gamma ) ) \\ 
&\hskip4mm + ( \text{terms from more split Cartan subgroups}),
\end{aligned}
\end{equation}
where   $p_j (u) \in \mathbb{Z}[u, u^{-1}]$. Note that  each root $\beta_k$ is simple for the preceding parameter 
$s_{\beta_{k-1}} \times \cdots \times (s_{\beta_1}\times \gamma)$.
 By \cite{Parkcity}, we can define the coherent continuation action of $W(\lambda)$ on 
$\bbZ[\calB]$, denoted $w\cdot \gamma$, with $w\in W(\lambda), \gamma\in \calB$, as follows. 

\smallskip

For $s_{\alpha}\in W(\lambda)$ 
with $\alpha\in \Pi(\lambda)$,
$$ 
  s_{\alpha}\cdot \gamma := -T_{s_{\alpha}} (\gamma)|_{u=1}, \ 
\text{with\,each\,term $\delta$ on\,the\,right\,side\,multiplied\,by $(-1) ^{ \ell (\gamma) -\ell  (\delta)}$,}
$$
where  $\ell$ is the length function defined on parameters (see Definition \ref{d:length}).

Therefore, from each step $T_{\alpha_j}$ in (\ref{hecke}), we may define
\begin{align*}
s_{\beta_j} \cdot \delta  &= -T_{  \beta_j } (\delta)| _{u=1}, \text{ if }  \alpha_ j  \text{ is real or imaginary for } \delta;  \\
s_{\beta_j} \cdot \delta  &=T_{  \beta_j } (\delta)|_{u=1}, \text{ if }  \alpha_j \text{ is  complex for } \delta.
\end{align*}

\begin{notation}
 Let $\alpha\in\Pi (\lambda)$.  Retain the notation in Lemma \ref{l:chamber}. 
Denote $m({\gamma}, s_{\alpha})$  the number of occurrences of imaginary roots in $ \{ \beta_j , 1\le j \le 2l+1  \}$ with respect
to the parameters  $s_{\beta_{j-1}} \times \cdots \times (s_{\beta_1}\times \gamma)$. 
\end{notation}

An easy calculation shows that 
\begin{equation} \label{key}
s_{\alpha}\cdot \gamma = (-1)^{m (\gamma,s_{\alpha}) } s_{\alpha}\times \gamma +\text{(terms from more split Cartan subgroups)}.
\end{equation}

Now fix a block $\calB _{\la,\chi}$ of regular characters of $\tu{G}$, then $W(\lambda)$ acts on
$\bbZ[\calB _{\lambda,\chi}]$ by the  coherent continuation action, since  $w\times \gamma \in \calB _{\la,\chi}$ for all $w\in W(\lambda), \gamma\in \calB$. Due
to the reason stated in the beginning of the section, the goal is to compute 
$[sgn_{W(\lambda)}   :  \bbZ[\calB_{\lambda,\chi}]]$, the multiplicity of the sign representation in $\bbZ[\calB _{\lambda,\chi}]$ when considered as $W(\lambda)$-representations.

\smallskip

Notice that two $\lambda$-regular characters $\gamma _i =(\widetilde{H_i}, \Gamma _i, \overline{\gamma_i})$ and $\gamma _j =(\widetilde{H_j}, \Gamma _j, \overline{\gamma_j})$ from $\mathcal{B}$ are in the same cross action orbit if and only if $\widetilde{H_i} = \widetilde{H_j}$. We enumerate the Cartan subgroups of $\widetilde{G}$ as $\{
\widetilde{H_1}, \cdots , \widetilde{H_m}\}$, and pick a  regular
character $\gamma _j$ specified by $\widetilde{H_j}$, then $\{\gamma
_1, \cdots, \gamma _m \}$ is a set of representatives of the cross
action orbits of $W(\lambda)$ on  $\mathbb{Z}[\mathcal{B_{\lambda,\chi}}]$.

\medskip

Let $W_{\gamma _j} = \{ w \in W(\lambda) \thinspace | \thinspace
w\times \gamma_j =\gamma _j  \}$ be the cross stabilizer of $\gamma _j
$ in $W(\lambda)$. 
Then we have the following proposition.

\begin{prop}  \label{p:coh-decomp}
As a $W(\lambda)$-representation, $\mathbb{Z}[\mathcal{B _{\lambda,\chi}}] \simeq \bigoplus _j  Ind ^{W(\lambda)}
_{W_{\gamma _j}   }  (\epsilon _j)$, where $\epsilon _j$ is a
one-dimensional representation of $W_{\gamma _j}$ such that for $w\in
W_{\gamma_j}$, $w \cdot \gamma _j = \epsilon _j (w) \gamma _j +$ other
terms from more split Cartan subgroups. 
\begin{proof}
This can be easily proved by the formulas given in \cite{RT} and (\ref{key}).
\end{proof}
\end{prop} 

By Proposition \ref{p:coh-decomp} and Frobenius reciprocity, the multiplicity of $sgn_{W(\lambda)}$ in $\mathbb{Z} [\calB_{\la,\chi}]$ 
is  $[sgn_{W(\lambda)}:   \mathbb{Z}[\calB_{\la,\chi}]] =[sgn_{W(\lambda)} | _{W_{\gamma_j}} : \epsilon _j   ]$, which is equal 
to 0 or 1, since $sgn_{W(\lambda)} | _{W_{\gamma_j}} $ is one-dimensional. This means that we have reduced our goal to count 
the number of $\gamma _j$'s making $[sgn_{W(\lambda)} | _{W_{\gamma_j}} : \epsilon _j   ] = 1$. Equivalently, we calculate 
the number of $\gamma_j$ such that  
\begin{equation}  \label{cond-star}
sgn_{W(\la)} |_{W_{\gamma_j}} = \ep_j.
\end{equation}
Due to Proposition \ref{p:coh-decomp} and (\ref{cond-star}), we have to analyze $W_{\gamma_j}$ and $\ep_j$ for each $\gamma_j$.

\medskip

Now we specialize to the case that $\tu G$ is of type $E$. 
 In this case, $\la=\rho/2$. We summarized 
$\rho/2$, $\Delta (\rho/2),$ and $\Pi(\rho/2)$ for each type in \ref{s:char-E}. 

We study the coherent continuation representation of $W(\rho/2)$ on $\bbZ [\calB _{\rho/2,\chi}]$ with a given genuine central character $\chi$ of $\tu G$. 

Recall the notation $W^r(\lambda), W^i (\lambda), W^C(\lambda)^{\theta}$ from Section 7 in \cite{AT}. For $\gamma\in \calB _{\rho/2, \chi}$, let 
$\Delta ^r$ and $\Delta ^i$ be the real and imaginary roots, and $\Delta ^r(\rho/2) = \Delta ^r \cap \Delta (\rho/2), \Delta ^i(\rho/2) = \Delta ^i \cap \Delta (\rho/2).$ Let 
$$
\Delta ^C (\rho/2) = \{\alpha \in \Delta (\rho/2) \mid \langle  \alpha ,\rho_r ^{\vee} \rangle =  \langle  \alpha ,\rho_i^{\vee} \rangle =0  \},
$$
where $\displaystyle \rho_r = \frac{1}{2} \sum _{\alpha \in \Delta ^+\cap \Delta ^r(\rho/2)} \alpha, \  \displaystyle \rho_i = \frac{1}{2} \sum _{\alpha \in \Delta ^+\cap \Delta ^i(\rho/2)} \alpha.$ 
\medskip

Let $W^r(\rho/2) =W(\Delta ^r (\rho/2)),
W^i(\rho/2)=W(\Delta ^i (\rho/2))$, and let $W^C(\rho/2) ^\theta$ be the fixed points of $\theta $ acting on $W^C(\rho/2) = W(\Delta ^C(\rho/2))$. 
 According to Proposition 7.14 in \cite{AT}, we have 

\begin{equation}
W_{\gamma}= [W^i (\rho/2)\times W^r (\rho/2)] \rtimes W^C(\rho/2) ^\theta
\end{equation}

We pick a parameter $\gamma _{r,m,s}$ corresponding to each Cartan subgroup $\tu H ^{r,m,s}$ to be the representative of the cross action orbit  as follows. 
Here we use the notation introduced in \ref{s:char-E}.
\begin{equation}
\begin{aligned}
E_6 :\  &  \gamma_{2,2,0} =\gamma (\{ 1, 2, 3, 4\})\\ 
      &   \gamma_{0,3,0} = \gamma (\{ 1, 2\}, \{ -1, -2\}, \{ 3, 4 \} )\\
     &  \gamma_{0,2,2} =  \gamma (\{ 1, 2\}, \{ -1, -2\} ) \\
           & \gamma_{0,1,4} = \gamma (\{ 1, 2\} ) \\
            &  \gamma_{0,0,6} =\gamma (\emptyset)
       \end{aligned}
\end{equation}

\begin{equation}
\begin{aligned}
E_7 :\   &  \gamma_ {7,0,0}  =\gamma (\{ 1, 2, 3, 4, 5, 6\}, \{7, 8\})    \\    
       &   \gamma _{5,1,0}=\gamma (\{ 1, 2, 3, 4\}, \{ 5, 6\}, \{7, 8\})  \\
        &   \gamma _{3,2,0}= \gamma (\{ 1, 2, 3, 4\},  \{7, 8\})   \\   
              &   \gamma_ {1,3,0}= \gamma (\{ 1, 2\},  \{ -1, -2\}, \{ 3, 4\},  \{7, 8\})   \\ 
              & \gamma_{2,2,1}= \gamma (\{ 1, 2, 3, 4\})    \\      
              &   \gamma_ {1,2,2}  =\gamma (\{ 1, 2\},  \{ -1, -2\},  \{7, 8\})  \\  
      &  \gamma_{0,3,1} =\gamma (\{ 1, 2\},  \{ -1, -2\},  \{3,4 \})    \\ 
 &  \gamma_{0,2,3}  =\gamma (\{ 1, 2\},  \{ -1, -2\})       \\  
 & \gamma_{0,1,5} =\gamma (\{ 1, 2\})     \\  
    &  \gamma_{0,0,7} =\gamma (\emptyset) \\  
                 \end{aligned}
\end{equation}
            \begin{equation}
\begin{aligned}
  E_8 :\  &    \gamma_{8, 0,0}=  \gamma (\{ 1, 2, 3, 4, 5, 6, 7, 8\}) \\     
            &     \gamma_{6,1,0} =  \gamma (\{ 1, 2, 3, 4, 5, 6\} , \{7, 8\}) \\     
             &  \gamma_{4,2,0} =  \gamma (\{ 1, 2, 3, 4\} , \{5, 6\} , \{7, 8\}) \\   
              &  \gamma_{2,3,0}= \gamma (\{ 1, 2, 3, 4\} , \{5, 6\} )\\  
               & \gamma_{0,4,0}=    \gamma (\{ 1, 2\}, \{-1,-2\}, \{3, 4\} , \{5, 6\} ) \\      
                 &  \gamma_{2,2,2}=\gamma (\{ 1, 2\}, \{-1,-2\}, \{3, 4\} , \{-3,-4\} ) \\  
                   &  \gamma_{0,3,2}=\gamma (\{ 1, 2\}, \{-1,-2\}, \{3, 4\} )   \\ 
                      &  \gamma_{0,2,4}= \gamma (\{ 1, 2\}, \{-1,-2\} )  \\  
                      &  \gamma_{0,1,6}= \gamma (\{ 1, 2\}) \\  
                       &  \gamma_{0,0,8}=\gamma (\emptyset) 
\end{aligned}
\end{equation}

In the following contents,  the notation for roots and Weyl group elements are from 
Notation \ref{n:E}. 

\begin{lemma}\label{l:ruleout1}
Suppose $\gamma$ is one of the following parameters:
\begin{align*}
E_6: \  &    \gamma_{0,2,2}, \gamma_{0,1,4}, \gamma_{0,0,6}\\
E_7: \ & \gamma_{0,0,7},\gamma_{0,1,5}, \gamma_{0,2,3},\gamma_{1,2,2}\\
E_8:\ & \gamma_{0,0,8}, \gamma_{0,1,6}, \gamma_{0,2,4}, \gamma_{0,3,2}, \gamma_{2,2,2}.
\end{align*}
Then $W^r (\rho/2) \neq 1$. More precisely, 
\begin{itemize}
\item[(a)] if $\gamma$ is one of the above parameters of type $E_6$ and $E_7$, $w=s_{35}\in W^r (\rho/2)$;
\item[(b)] if $\gamma$ is one of the above parameters of type $E_8$, $w=s_{ 57}\in W^r (\rho/2)$;
\end{itemize}
For $w$ chosen above, we have that $\epsilon _{\gamma}(w)=1$.
\begin{proof}
For the case of type $E_6$ and $E_7$, $\theta_\gamma(-\alpha_3+\alpha_5 ) =-(-\alpha_3+\alpha_5 )$ for each $\gamma$ 
listed, so $s_{35}\in W^r(\rho/2).$ Decompose $s_{35}$ as a product of simple reflections: $s_{35}=s_{45}s_{34}s_{45}$. Furthermore, we decompose
$s_{35}$ into a product in  (\ref{e:simple}) with $\beta_1=-e_4+e_5, \beta_2=-e_3+e_5, \beta_3 = -e_3+e_4$. Using 
$\theta_\gamma$ for each listed $\gamma$, it can be readily verified that $\theta_\gamma(\beta_i)=-\beta_i$ and hence each $\beta_i$ is a real
root for the parameter $s_{\beta_{i-1}} \times ( s_{\beta_{i-2} } \times \cdots (s_{\beta_1}\times \gamma))$. Consequently, $m(\gamma, s_{35} )=0$ 
and hence $w_{\gamma} (s_{35})=1$.

 Similarly for the case of $E_8$: $s_{57} = s_{67}s_{56} s_{67}$, and this is equal
to the decomposition in (\ref{e:simple}) with $\beta_1=-e_6+e_7, \beta_2=-e_5+e_7, \beta_3=-e_5+e_6$. The rest argument is the same as the case of 
$E_6$ and $E_7$. 

\end{proof}
\end{lemma}


\begin{lemma}\label{l:ruleout2}
Suppose $\gamma$ is one of the following parameters:
\begin{align*}
E_6: \  &    \gamma_{0,3,0} \\
E_7: \ & \gamma_{5,1,0}, \gamma_{3,2,0}, \gamma_{1,3,0}, \gamma_{0,3,1}  \\
E_8:\ & \gamma_{6,1,0}, \gamma_{4,2,0}, \gamma_{2,3,0}, \gamma_{0,4,0}.
\end{align*}
Then there exists an element $w\in W^C(\rho/2)^\theta $ such that $\epsilon_\gamma (w)=-1$. More precisely,
$w$ can be chosen as follows:
  \begin{align*}
E_6: \  &    \gamma_{0,3,0}   & w= s_{13} s_{\Br{13}} s_{24} s_{\ovl{24}} \\
E_7: \ & \gamma_{5,1,0}, \gamma_{1,3,0}, \gamma_{0,3,1}  &w= s_{13} s_{\Br{13}} s_{24} s_{\ovl{24}}  \\
 &   \gamma_{3,2,0}   &   { w=  r_{4567}r_{123567}    r_{1247} r_{37}  } \\
E_8:\ & \gamma_{6,1,0}, \gamma_{4,2,0}, \gamma_{2,3,0}, \gamma_{0,4,0} & w= s_{13} s_{\Br{13}} s_{24} s_{\ovl{24}}.
\end{align*}

\begin{proof}
Let $\gamma$ be one of the parameters listed above. 

First suppose $\gamma$ is not $\gamma_{3,2,0}$ of type $E_7.$  Let $w=s_{13}s_{\Br{13}}s_{24}s_{\Br{24}}$. Note that 
$\theta_{\gamma}=s_{12} s_{\Br{12}}s_{34} w'$ for some $w'\in W$ with $ww'=w'w$. So to show $w\in W^C(\rho/2)^\theta$
it suffices to show that $s_{12} s_{\Br{12}}s_{34} w=w s_{12} s_{\Br{12}}s_{34}$. This equality holds since both the right hand
side and left hand side maps $(x_1,x_2,x_3,x_4)$ to $(x_1,x_2,x_4,x_3)$ for any sequence of numbers $(x_1,x_2,x_3,x_4)$. 
Now we decompose $w=s_{13}s_{\Br{13}}s_{24}s_{\Br{24}}$ into a product of simple roots:
\begin{align*}
w &=  s_{13}s_{\Br{13}}s_{24}s_{\Br{24}}\\
    &= (s_{23} s_{12} s_{23}) (s_{23} s_{\Br{12}} s_{23}) (s_{34} s_{23} s_{34}) (s_{34} s_{12} s_{23} s_{\Br{12}} s_{23} s_{12} s_{34})\\
    & =  s_{23} s_{12}  s_{\Br{12}} s_{23} s_{34} s_{23}  s_{12} s_{23} s_{\Br{12}} s_{23} s_{12} s_{34} 
\end{align*}
Write this to be a product  $w=s_{\beta_{12}} \cdots w_{\beta_1}$ as in (\ref{e:simple}), where  the $\beta_i$'s are in turn 
\begin{eqnarray*}
-e_3+e_4 \ (im), \ -e_1+e_2 \ (im), \  -e_1+e_4 \ (cx) \ ,  e_2+e_4 \ (cx),  \ e_1+e_2 \ (im), \  e_1+e_4\ (cx)\\
e_2-e_4 \ (cx), \ e_3+e_4\ (real), \ e_2+e_3\ (cx), \ e_1+e_3\ (cx), \ -e_1+e_3\ (cx), \ -e_2+e_3\  (cx).
\end{eqnarray*}
The word inside parentheses after every $\beta_i$ indicates the root type with respect to the preceding parameter. 
Therefore, $m(\gamma, w)=3$, and hence $\epsilon _{\gamma} (w)=-1$. 
\medskip

Now suppose that $\gamma=\gamma_{3,2,0}$ with type $E_7$.  Take $w = r_{4567} r_{123567} r_{1247} r_{37}$.
 Since $\theta= \theta_\gamma = (-1) s_{12}s_{\Br{12}} s_{34} s_{\Br{34}} s_{78}$, 
 \begin{eqnarray*}
 \theta w \theta^{-1} &=& s_{\theta \beta_{4567}}   s_{\theta \beta_{123567}}  s_{\theta \beta_{1247}}  s_{\theta \beta_{37}} \\
&=& r_{47} r_{1237} r_{124567} r_{3567}.
\end{eqnarray*}

 To show $w\in W^C(\rho/2)^\theta$, we have to show 
$r_{4567} r_{123567} r_{1247} r_{37} =r_{47} r_{1237} r_{124567} r_{3567}$. This is a straightforward calculation by checking that 
 these two Weyl group elements agree on the simple roots of $E_7$. 

Now we decompose $w =  r_{4567} r_{123567} r_{1247} r_{37}$ into a product of simple roots:
\begin{align*}
r_{4567} &=  s_{12}s_{23} s_{\Br{12}} s_{23} s_{12}  r_{234567}    s_{12}s_{23} s_{\Br{12}} s_{23} s_{12}  \\
   r_{123567} &=s_{34}s_{23} s_{12}  r_{234567}    s_{12}s_{23}s_{34}  \\
    r_{1247}& =  s_{ \Br{56}} s_{23}  s_{12}  r_{234567} s_{12} s_{23} s_{\Br{56}} \\
    r_{37}&= s_{\Br{56}} s_{\Br{24}} r_{234567} s_{\Br{24}} s_{\Br{56}},
\end{align*}
where $s_{\Br{56}} = s_{45} s_{34} s_{23} s_{12} s_{56} s_{45} s_{34} s_{23}s_{\Br{12}} s_{23} s_{34} s_{45}s_{56} s_{12} s_{23} s_{34} s_{45}$ and 
$s_{\Br{24}} = s_{12} s_{34} s_{23} s_{\Br{12}} s_{23} s_{34} s_{12}$.
Write this to be a product  $w=s_{\beta_{72}} \cdots s_{\beta_1}$ as in (\ref{e:simple}), where  the $\beta_i$'s are in turn  (up to a sign)
\begin{eqnarray*}
&-e_4+e_5 \ (cx), \ -e_3+e_5 \ (cx), \  -e_2+e_5 \ (cx) \ ,  -e_1+e_5 \ (cx),  \ -e_4+e_6 \ (cx), \  -e_3+e_6 \ (cx)\\
&-e_2+e_6\ (cx), \ -e_1+e_6 \ (cx), \  e_5+e_6 \ (real) \ ,  e_1+e_5 \ (cx),  \ e_2+e_5\ (cx), \  e_3+e_5\ (cx)\\
&e_4 +e_5 \ (cx), \ e_1+e_6 \ (cx), \  e_2+e_6 \ (cx) \ ,  e_3+e_6 \ (cx),  \ e_4+e_6 \ (cx), \  -e_1+e_2\ (im)\\
&-e_3+e_4 \ (im), \ -e_1+e_4 \ (im), \  e_2+e_4 \ (im) \ ,  e_1+e_2 \ (im),  \ e_2+e_3\ (im), \  e_1+e_4\ (im)\\
&\beta_{37} \ (cx), \ \beta_{1347}\ (cx), \  e_2+e_3 \ (im) \ ,  \beta_{1237}\ (cx),  \ \beta_{2347}\ (cx), \  -e_1+e_4\ (im)\\
&\beta_{47} \ (cx), \  -e_1+e_2\ (im), \  \beta_{27} \ (cx) \ ,  \beta_{17}\ (cx),  \ \beta_{1247}\ (cx), \  e_2+e_4\ (im)\\
&e_1+e_4 \ (im), \  \beta_{3467}\ (cx), \  \beta_{123467} \ (cx) \ ,  \beta_{2367}\ (cx),  \ \beta_{1367}\ (cx), \  \beta_{3457}\ (cx)\\
 &\beta_{123457}\ (cx), \  \beta_{2357}\ (cx), \  \beta_{1357} \ (cx) \ ,  -e_7+e_8\ (im),  \ \beta_{2467}\ (cx), \  \beta_{1467}\ (cx)\\
 &\beta_{67}\ (cx), \  \beta_{1267}\ (cx), \  \beta_{2457} \ (cx) \ ,  \beta_{1457}\ (cx),  \ \beta_{57}\ (cx), \  \beta_{1257}\ (cx)\\
&e_1+e_2 \ (im), \ -e_2+e_4 \ (im), \  -e_1+e_4 \ (im) \ ,  \beta_{123567} \ (cx),  \ \beta_{234567}\ (cx), \  \beta_{134567}\ (cx)\\
&\beta_{3567} \ (cx), \ -e_1+e_2 \ (im),\  e_2+e_4 \ (im) \ ,  \beta_{1567} \ (cx),  \ \beta_{124567}\ (cx), \  \beta_{2567}\ (cx)\\
&\beta_{4567} \ (cx), \ -e_2+e_4 \ (im),\  e_1+e_2 \ (im) \ ,  -e_1+e_4 \ (im),  \  e_2+e_4\ (im), \  e_1-e_2\ (im).
\end{eqnarray*}
Here \emph{im} means the root is imaginary, and \emph{cx} means the root is complex. It turns out that 
there are 23 imaginary roots showing up in this list.
Therefore, $m(\gamma, w)=23$, and hence $\epsilon _{\gamma} (w)=-1$. 
\end{proof}

\end{lemma} 
\begin{remark} \label{r:ruleout}
For $w\in W^r(\rho/2)$ in Lemma \ref{l:ruleout1}, $sgn (w) =-1$, and hence the corresponding parameter $\gamma$ fails to satisfy (\ref{cond-star}). 
For $w\in W^C(\rho/2)^\theta$ in Lemma \ref{l:ruleout2}, $sgn (w) =1$, and hence the corresponding parameter $\gamma$ fails to satisfy (\ref{cond-star}). 

\end{remark}
\begin{theorem} \label{t:Ecounting}
Suppose $\tu G$ is the split group of type $E$. Then $| \prod _{\rho/2} ^s (\tu G)| = 1, 4, 1$, if $\tu G$ is of type $E_6, E_7, E_8,$ respectively.  Accordingly, 
$ \prod _{\rho/2}  ^s (\tu G) = \prod _{R_D} (\tu G).$
\begin{proof}
Due to the observation in Remark \ref{r:ruleout}, the parameters in Lemma \ref{l:ruleout1} and \ref{l:ruleout2} fail to satisfy (\ref{cond-star}). Therefore, the parameters for each type  (with a fixed central character) which have not been
ruled out are
\begin{align*}
E_6: \  &    \gamma_{2,2,0} \\
E_7: \ & \gamma_{7,0,0}, \gamma_{2,2,1} \\
E_8:\ & \gamma_{8,0,0}.
\end{align*}
This means that $|\prod _{\rho/2} ^s (\tu G)|\le 1, 4, 1$ for $E_6, E_7, E_8$, respectively, 
since the number of genuine central characters is 1 for $E_6$, 2 for $E_7$, and 1 for $E_8$.
In \cite{Ts}, it is shown that  $\prod _{R_D} (\tu G) \subseteq  \prod _{\rho/2} ^s (\tu G)$ and $|\prod _{R_D} (\tu G)| = 1, 4, 1$ for $E_6, E_7, E_8$, respectively. Therefore,
$| \prod _{\rho/2} ^s (\tu G)| = 1, 4, 1$ for $E_6, E_7, E_8,$ respectively, and  $ \prod _{\rho/2} ^s (\tu G) = \prod _{R_D} (\tu G).$
\end{proof}
\end{theorem}

\section{Proof of the main Theorem}\label{s:main}

Let $\gamma_0$ be the regular character corresponding to the trivial representation $\bbC$ and let $I_G(\gamma_0)$ be the standard module parametrized by $\gamma_0$. 
Write $\bbC$ as a sum of standard modules, say, 
\begin{equation} \label{e:Z}
\bbC = \sum\limits _\gamma M(\gamma ,\gamma_0) I_G(\gamma),
\end{equation} 
with $\gamma$ taken over holomorphic characters and $M(\gamma, \gamma_0)=(-1)^{\ell (\gamma)-\ell(\gamma_0)}$ (see \cite[Section 7]{A1}, for example). 

Equation (\ref{e:Z}) is the Zuckerman character formula for the trivial representation. We take Lift$_G ^{\tu G}$ on (\ref{e:Z}). After lifting, some terms on the right vanish:
\begin{equation} \label{e:liftZ}
\text{Lift}_G ^{\tu G} (\bbC) = \sum \limits_{\text{some} \  \tu \gamma} c_{\tu \gamma} I_{\tu G} (\tu\gamma).
\end{equation}

To prove the main theorem, it amounts to showing that $J(\tu \delta)$ shows up on the right side of (\ref{e:liftZ}) for $J(\tu \delta) \in \prod _{R_D}(\tu G)$. 
\medskip

Let $Sh\in\cal{SH}$ and let $\gamma_{Sh}$ be the corresponding regular character. Every representation in $\prod _{R_D}(\tu G)$ is of the form 
$J(c_S(\gamma_{Sh}))$ for some $S\in R_D$. The following Lemma characterizes the standard modules $I_{\tu G} (\tu\delta)$ which contains 
$J(c_S(\gamma_{Sh}))$ as a composition factor. 

Recall that $P_{\gamma,\delta} (u)$ is a Kazhdan-Lusztig polynomial and $M(\gamma,\delta) = (-1)^{\ell(\delta) -\ell (\gamma)}  P_{\gamma,\delta} (1)$. 
\begin{lemma}\label{l:keyA}
Using the notation in Section \ref{s:typeA}, let
\begin{itemize}
\item[(a)]  $\gamma^*=\gamma(\{1,2\}, \dots, \{2p-1,2p\})$ if  $\tu G=\tu{SL}(2p,\bbR)$;
\item[(b)]  $\gamma^* =\gamma(\{1,2\}, \dots, \{2p-1,2p\}),  \ \gamma(\{1,2\}, \dots,  \{ 2p-3, 2p-2     \} ,\{-(2p-1),-2p\}),$ or $\gamma(\{ 2p-1, 2p\} , \{-(2p-1), -2p\} )$
if $\tu G=\tu{Spin}(2p,2p)$;
\item[(c)]  $\gamma^* = \gamma(\{ 2p-1, 2p\} , \{-(2p-1), -2p\} )$ if $\tu G=\tu{Spin}(2p+1,2p+1)$;
\item[(d)] $\gamma^* = \gamma (\{-1,-2\} ,\{3, 4\} ,\{5,6\}) $ if $\tu G=\tu{E_7}$\emph{(split)}.
\end{itemize}
 Then  $P_{\gamma^*, \delta} (1)=
\begin{cases}
1 & \text{ if } \delta\subseteq  \gamma^*, \\
0 & \text{ if } \delta \not\subseteq \gamma^*.
\end{cases}
$
Consequently, 
$M (\gamma^*,\delta)=
\begin{cases}
(-1)^{\ell (\delta)-\ell(\gamma^*)} & \text{ if } \delta\subseteq  \gamma^*, \\
0 & \text{ if } \delta \not\subseteq \gamma^*.
\end{cases}
$
\begin{proof}


We prove the assertion for (a). The argument for the rest cases is similar.

Let $\delta \subsetneq \gamma^*$. Suppose that $\{2i_1-1, 2i_1 \} , \{2i_2-1, 2i_2 \}, \dots, \{2i_k-1, 2i_k \} $ are the pairs missing from $\delta$. This means that 
$\alpha_1= e_{2i_1-1}-e_{2i_1}, \dots, \alpha_k=e_{2i_k-1}-e_{2i_k}$ are imaginary for $\gamma^*$ and are real for $\delta$. Write $s_1,\dots, s_k$ to be the simple reflections corresponding 
to $\alpha_1,\dots,\alpha_k$, respectively. Then by  by Case III (iv) of \cite[Proposition 7.10]{RT}
\begin{eqnarray*}
P_{\gamma^*, \delta } (1) &=&  P_{s_1\times \gamma^*, (s_1\times \delta)_{\alpha_1} } (1)  = P_{ s_2\times (s_1\times \gamma^*), (s_2\times(s_1\times \delta)_{\alpha_1} )_{\alpha_2}} (1) \\
&=&\cdots \cdots = P_{ s_k \times (\cdots \times s_2\times (s_1\times \gamma^*)), (s_k\times \cdots (s_2\times(s_1\times \delta)_{\alpha_1} )_{\alpha_2} \cdots )_{\alpha_k} } (1) \\
&=&  P_{ s_k \times (\cdots \times s_2\times (s_1\times \gamma^*)),s_k \times (\cdots \times s_2\times (s_1\times \gamma^*)) } (1) =1
\end{eqnarray*}

Now suppose that $\delta \not\subseteq \gamma^*$ with $\ell (\delta)>\ell (\gamma^*)$. 
Suppose that $\{i_1, j_1 \} , \{ i_2,  j_2 \}, \dots, \{i_k, j_k\} $ are the pairs missing from $\delta$, where the $i$'s are in $\{1,3,\dots,2p-1\}$, $j$'s are in $\{2,4,\dots,2p\}$.
Note that $k\ge 1$; otherwise, $\ell(\delta)<\ell (\gamma^*)$.
Then
$\alpha_1= e_{i_1}-e_{j_1}, \dots, \alpha_k=e_{i_k}-e_{j_k}$ are complex for $\gamma^*$ and  are real for $\delta$. 
Each $\alpha_l$ is not necessarily simple in $\Delta^+$. We move to another chamber if needed, say, $w\Delta^+$ with $w=t_m\cdots t_2t_1$, where $t_{i+1}$ is simple 
in $t_i(\Delta^+)$, such that $\alpha_1,\dots, \alpha_k$ are simple in $w\Delta^+$. 
Define 
\begin{align*}
\gamma^\# = t_m\times (\cdots \times t_2\times( t_1\times \gamma^*)) , \text{ and } \\
\delta^\# = t_m\times (\cdots \times t_2\times( t_1\times \delta)).
\end{align*}
In this way, $\alpha_1,\dots,\alpha_k$ are again complex for $\gamma^\#$ and are real for $\delta ^\#$.
Write $s_1,\dots, s_k$ to be the simple reflections corresponding 
to $\alpha_1,\dots,\alpha_k$, respectively. Then  by Case III (i) and (ii) of \cite[Proposition 7.10]{RT},
\begin{eqnarray*}
P_{\gamma^*, \delta } (1) &=&  P_{\gamma^\#, \delta^\#} (1)
= P_{s_1\times \gamma^\#, (s_1\times \delta^\#)_{\alpha_1} } (1)  = P_{ s_2\times (s_1\times \gamma^\#), (s_2\times(s_1\times \delta^\#)_{\alpha_1} )_{\alpha_2}} (1) \\
&=&\cdots \cdots = P_{ s_k \times (\cdots \times s_2\times (s_1\times \gamma^\#)), (s_k\times \cdots (s_2\times(s_1\times \delta ^\#)_{\alpha_1} )_{\alpha_2} \cdots )_{\alpha_k} } (1) \\
&=&0
\end{eqnarray*}
The last equality is true since $$\ell (s_k\times \cdots (s_2\times(s_1\times \delta^\#)_{\alpha_1} )_{\alpha_2} \cdots )_{\alpha_k} ) < \ell (s_k \times (\cdots \times s_2\times (s_1\times \gamma^\#) ).$$
\end{proof}
\end{lemma}

\begin{lemma}\label{l:inverseKL}
Fix  $S\in R_D$ and fix $Sh\in \cal{SH}$ with central character $\chi$. Then $J(c_S(\gamma_{Sh}))$ is a composition factor 
of $I_{\tu G} (\tu \gamma)$ if and only if $I_{\tu G} (\tu \gamma)$ is of the form $I_{\tu G} (c_{S'} (\gamma_{Sh}))$,  where $S'\subseteq S$. 
Moreover, the multiplicity of $J(c_S(\gamma_{Sh}))$ in $I_{\tu G} (c_{S'} (\gamma_{Sh}))$ is 1, i.e. $$m(c_S (\gamma_{Sh}) , c_{S'} (\gamma_{Sh})) = 1.$$
\begin{proof}

It suffices to show that for $\gamma^*$ listed in Lemma \ref{l:keyA}, we have that 
\begin{itemize}
\item[(i)] $m(\gamma^* , \delta ^*)=1$ if $\delta ^*\subseteq \gamma^*$;
\item[(ii)]  $m(\gamma^* , \delta^*)=0$ if $\delta^* \not\subseteq \gamma^*$.
\end{itemize}

When $\tu G=\tu{SL}(2n,\bbR)$, 
$\gamma^*=\gamma (\{1,2\} , \{3,4\}, \dots, \{2p-1,2p\}).$

We prove (i) by induction on number of pairs of $\{i, i+1\}$ (with $i$ odd) occurring in $\gamma^*$. 

Suppose $\gamma^*=\gamma(\{1,2\})$. 
By Case III (iv) of Proposition 7.10  in \cite{RT}, $P_{\gamma^* ,\gamma (\emptyset)} (1)=1$, and hence $M(\gamma^* , \gamma (\emptyset)) =-1$. 
Then
\begin{eqnarray*}
J(\{1,2\}) &=& I(\{1,2\}),\\
J(\emptyset)&=& I(\emptyset) -I(\{1,2\}).
\end{eqnarray*} 
Consequently, $I(\emptyset) = J(\emptyset) + J(\{1,2\})$, meaning that $m( \gamma^* , \gamma (\emptyset)) =1$ and $m( \gamma^* , \gamma^*) =1$.

\medskip
Now for any $\gamma^\diamond = \gamma(  \underbrace{\{1,2\} , \dots, \{2p-1,2p\} } _{\# \text{ of pairs } <p})$ with $\gamma^\diamond \subsetneq \gamma ^*$, 
suppose that  $m(\gamma^\diamond, \delta)=1$ if $\delta\subseteq \gamma^\diamond$. 

Fix $\delta^*$ with $\delta^*\subsetneq \gamma ^*.$  Let $\ell (\delta^*)=\ell_s-q$ with $q<p$. Since $M(\gamma,\delta) m(\gamma,\delta)=Id$, we have 
\begin{equation}
\sum_{\delta} M(\gamma^*,\delta) m (\delta, \delta^*) =0.
\end{equation}
This can be rewritten as 
\begin{equation} \label{e:Mm}
0=M(\gamma^*,\gamma^*) m(\gamma^*,\delta^*) + \sum_{\substack{\delta\\ \ell(\gamma^*) <\delta <\ell(\delta ^*) }} M(\gamma^*,\delta) m(\delta,\delta^*) + M(\gamma^*,\delta^*) m(\delta^*,\delta^*).
\end{equation}
Therefore,
\begin{equation}
m(\gamma^*,\delta^*)= -\left ( \sum _{\substack{ \delta \\ \delta^* \subset \delta \subset \gamma^*}}  M(\gamma^*,\delta) m (\delta, \delta^*)   \right ) - (-1) ^{\ell (\gamma^*)-\ell (\delta ^*)}.
\end{equation}
By induction hypothesis, $\{ \delta\mid m(\delta, \delta ^*) =1\} = \{ \delta \mid \delta^*\subseteq \delta \}$, along with Lemma \ref{l:keyA}, we have
\begin{equation}
m(\gamma^*,\delta ^*)= -\left ( \sum _{i=1} ^{p-q-1} (-1)^i {p-q \choose i} \right ) - (-1)^{p-q} =1 
\end{equation}
since $-\left ( \sum _{i=1} ^{p-q-1} (-1)^i {p-q \choose i} \right ) =\begin{cases} 0 & \text{ if $p-q$ is odd,} \\ 2 & \text{ if $p-q$ is even. } \end{cases}$

Next we prove (ii) by induction on $p$. For any $\gamma^\diamond =
\gamma(  \underbrace{\{1,2\} , \dots, \{2p-1,2p\} } _{\# \text{ of pairs } <p})$ with $\gamma^\diamond \subsetneq \gamma ^*$,
suppose that $m(\gamma^\diamond,\delta) =0$ if $\delta \not\subseteq \gamma^\diamond$.

Now fix $\delta^*$ with $\delta^*\not\subseteq \gamma ^*$ and  $\ell (\delta^*)=\ell_s-q$ with $q<p$.
We have equation (\ref{e:Mm}) again, whereas in this case (\ref{e:Mm}) becomes
\begin{equation} \label{e:Mm2}
m(\gamma^*,\delta^*)= -\left ( \sum _{\substack{ \delta \\  \delta \subsetneq \gamma^*}}  M(\gamma^*,\delta) m (\delta, \delta^*)   \right )
\end{equation}

In the sum on the right side of (\ref{e:Mm2}), since $\delta\subset \gamma^*$, we have $\delta ^*\not\subseteq \delta$, and hence $m(\delta,\delta^*) =0$
by induction hypothesis for every $\delta$ occurring in (\ref{e:Mm2}). 
Consequently, $m(\gamma^*,\delta^*)=0$.

The argument for $\gamma^* = \gamma(\{1,2\}, \dots,  \{ 2p-3, 2p-2     \} ,\{\ep(2p-1),\ep 2p\})$, $\ep=\pm 1$, when $\tu G=\tu{Spin}(2p, 2p)$
is analogous to that of $\tu G=\tu{SL}(2n,\bbR)$. 

Now suppose $\gamma ^*=\gamma(\{2p-1, 2p+1\} , \{-(2p-1), -(2p+1)\})$ in type $D_{2p}$. Then $\ell (\gamma^*)=\ell (\gamma (\emptyset)) -2$. 
We just need to show (i) and (ii) for $\delta^*$ with $\ell (\delta^*)=\ell (\gamma (\emptyset)) -1$. Fix such $\delta ^* \subset \gamma^*$, that is, $\delta^*=\gamma(\{2p-1, 2p+1\})$ or 
$\gamma ( \{-(2p-1), -(2p+1)\})$. Then $M(\gamma^*,\delta^*)=-1$ (by the argument in the proof of Lemma \ref{l:keyA}). Therefore, (\ref{e:Mm}) becomes 
$$
0=M(\gamma^*,\gamma^*)m(\gamma^*,\delta^*) + M(\gamma^*,\delta^*) m(\delta^*,\delta^* )= m(\gamma^*,\delta^*)-1,
$$ 
and hence $m(\gamma^*, \delta^*) =1$. 

If $\delta^*\not\subset\gamma^*$ with $\ell (\delta^*)=\ell (\gamma (\emptyset)) -1$, we have $M(\gamma^*,\delta^*)=0$, and hence (\ref{e:Mm}) gives
$0=M(\gamma^*,\gamma^*)m(\gamma^*,\delta^*)$, so $m(\gamma^*,\delta^*)=0,$ as desired. 

The case of  $\gamma^*= \gamma (\{ -1, -2\}, \{3,4\} , \{5,6\})$ in type $E_7$ is similar.
 
\end{proof}

\end{lemma}

\begin{lemma} \label{l:trivKL}
Fix $S\in R_D$. Let $\gamma_0$ be the regular character for $\bbC$. Then $$M(c_{S'} (\gamma_0), \gamma_0) =(-1)^{\ell(\gamma_0)  - \ell(c_{S'} (\gamma_0))}$$ for
$S'\subseteq S.$ 

\begin{proof}
This can be easily checked by the Kazhdan-Lusztig-Vogan algorithm for linear groups (see \cite{Parkcity}). 
\end{proof}
\end{lemma}

\begin{theorem}
Let $\tu G$ be split of type $A, D, E$. Then $\prod _{R_D} (\tu G) \subseteq \text{\emph{Lift}}(\bbC).$
\begin{proof}
We need to show that the coefficient of each $J(\tu \delta)\in \prod_{R_D}(\tu G)$ in (\ref{e:liftZ}) is nonzero. By Lemma \ref{l:inverseKL}, for each $S\in R_D$, one 
needs to compute the coefficient of $I(c_{S'} (\gamma_{Sh}))$ for $S'\subseteq S$. 

It is clear that every $Sh$ is in Lift$(\bbC)$ since the only standard module containing $Sh$  in (\ref{e:liftZ}) is $I(\gamma_{Sh})$. Therefore, we just need to show the theorem 
when $|\prod_{R_D} (\tu G) | >1$. That is, we will compute the coefficient of each $J(c_S (\gamma_{Sh}))$ in  (\ref{e:liftZ}) for $\emptyset\neq S\in R_D$.

\medskip

By Theorem  \ref{Liftstable},  Lemma \ref{l:inverseKL}  and  Lemma \ref{l:trivKL}, for $S\in R_D,$ the coefficient of $J(c_S(\gamma_{Sh}))$ in Lift$(\bbC)$ is 
\begin{equation} \label{e:coeff}
 \sum_{S'\subseteq S} (-1)^{\ell(\gamma_0) - \ell (\gamma_{S'})} C(H_{S'}), 
\end{equation}
which is denoted $K_S$, 
where $\gamma_0$ is the regular character for $\bbC$, $\gamma_{S'} := c_{S'}(\gamma_0)$ is a regular character specified by $H_{S'}$. We compute the constants
$C(H_{S'})$ (defined in Theorem \ref{Liftstable})  and the coefficients $K_S$ case by case.  We use the labeling in \cite[Table 3]{Ts} for the roots $\alpha_j$ in $S'$ and $S$. 

\begin{itemize}
\item[(i)] \underline{Type $A_{n-1}$, $n=2p$}

When $S= \{ \alpha_1, \alpha_3,\cdots, \alpha_{2p-1}  \}$:
\medskip
\begin{align*}
   S'\hspace{2.4em} &&\hskip-2mm  \#\{ S' \}   \hspace{1em}  &&  H_{S'} \   \  \ &&  \  \ \hskip-2mm  C(H_{S'}) \\
   \emptyset \hspace{2.7em}   &&  \ \ 1 \hspace{1em} \ \  &&\hskip-2mm  (\bbR^\times)^{n-1}&&  \hskip-2mm  1 \quad \\
 \{\alpha_i\}   \hspace{2em}&&  \ \ {p\choose 1} \hspace{1em}  &&  \hskip-2mm  \quad (\mathbb{R}^{\times})  ^{n-3} \times\mathbb{C}^{\times} &&  \hskip-2mm  1 \quad \\
\vdots \ \  \quad\quad  &&  \ \ \vdots \hspace{1.7em}   &&  \hskip-2mm  \vdots \quad\quad  &&  \hskip-2mm  \vdots \quad \\
  \{\alpha_{i_1},\dots, \alpha_{i_k}\}     &&  \ \ {p\choose k} \hspace{1em}  &&  \hskip-2mm  \quad (\mathbb{R}^{\times})  ^{n-1-2k} \times(\mathbb{C}^{\times} )^k&&  \hskip-2mm  1 \quad\\
 \vdots \   \ \quad \quad&&  \ \ \vdots  \hspace{1.7em}  &&  \hskip-2mm  \vdots \quad\quad  &&  \hskip-2mm  \vdots \quad \\
 S-\{\alpha_i\}   \hspace{1em}  &&  {p\choose p-1}    &&  \hskip-2mm  \quad \mathbb{R}^{\times} \times(\mathbb{C}^{\times} )^{p-1}&&  \hskip-2mm  1 \quad\\
  S  \hspace{2.4em} &&  \ \ 1 \hspace{1.5em} &&  \hskip-2mm  \quad (\mathbb{C}^{\times} )^{p-1}\times S^1&&  \hskip-2mm  2 \quad
  \end{align*}

In this table, $i\in \{1, 3, \dots, 2p-1 \}$ and $\{i_1, \dots, i_k\}\subset \{ 1, 3, \dots, 2p-1\}$. Then from (\ref{e:coeff}), 
$$
K_S= \sum_{k=0} ^{p-1} (-1)^k {p \choose k} +(-1)^p\cdot 2=
\begin{cases}
1& \text{ if $p$ is even,}\\
-1 & \text{ if $p$  is odd}. 
\end{cases}
$$

\item[(ii)] \underline{Type $D_{n}$, $n=2p$}

When $S= \{ \alpha_{n-1}, \alpha_n  \}$:
\medskip
\begin{align*}
   S'\hspace{2.4em} &&\hskip-2mm  \#\{ S' \}   \hspace{1em}  &&  H_{S'} \   \  \ &&  \  \ \hskip-2mm  C(H_{S'}) \\
   \emptyset \hspace{2.7em}   &&  \ \ 1 \hspace{1em} \ \  &&\hskip-2mm  (\bbR^\times)^{n}&&  \hskip-2mm  1 \quad \\
 \{\alpha_i\}   \hspace{2em}&&  \ \ 2 \hspace{1.7em}  &&  \hskip-2mm  \quad (\mathbb{R}^{\times})  ^{n-2} \times\mathbb{C}^{\times} &&  \hskip-2mm  1 \quad \\
  S  \hspace{2.4em} &&  \ \ 1 \hspace{1.7em} &&  \hskip-2mm  \quad (\bbR^\times)^{n-3}\times \mathbb{C}^{\times} \times S^1&&  \hskip-2mm  2 \quad
  \end{align*}

In this table, $i\in \{ n-1,n \}$. Then from (\ref{e:coeff}), 
$$
K_S 
=1-2\cdot 1 +2=1.$$

When $S= \{ \alpha_1, \alpha_3,\cdots, \alpha_{2p-3}, \alpha_{2p-1}  \}$ or 
$\{ \alpha_1, \alpha_3,\cdots, \alpha_{2p-3}, \alpha_{2p}  \}$ :
\medskip
\begin{align*}
   S'\hspace{2.4em} &&\hskip-2mm  \#\{ S' \}   \hspace{1em}  &&  H_{S'} \   \  \ &&  \  \ \hskip-2mm  C(H_{S'}) \\
   \emptyset \hspace{2.7em}   &&  \ \ 1 \hspace{1em} \ \  &&\hskip-2mm  (\bbR^\times)^{n}&&  \hskip-2mm  1 \quad \\
 \{\alpha_i\}   \hspace{2em}&&  \ \ {p\choose 1} \hspace{1em}  &&  \hskip-2mm  \quad (\mathbb{R}^{\times})  ^{n-2} \times\mathbb{C}^{\times} &&  \hskip-2mm  1 \quad \\
\vdots \ \  \quad\quad  &&  \ \ \vdots \hspace{1.7em}   &&  \hskip-2mm  \vdots \quad\quad  &&  \hskip-2mm  \vdots \quad \\
  \{\alpha_{i_1},\dots, \alpha_{i_k}\}     &&  \ \ {p\choose k} \hspace{1em}  &&  \hskip-2mm  \quad (\mathbb{R}^{\times})  ^{n-2k} \times(\mathbb{C}^{\times} )^k&&  \hskip-2mm  1 \quad\\
 \vdots \   \ \quad \quad&&  \ \ \vdots  \hspace{1.7em}  &&  \hskip-2mm  \vdots \quad\quad  &&  \hskip-2mm  \vdots \quad \\
 S-\{\alpha_i\}   \hspace{1em}  &&  {p\choose p-1}    &&  \hskip-2mm  \quad (\mathbb{R}^{\times})^2 \times(\mathbb{C}^{\times} )^{p-1}&&  \hskip-2mm  1 \quad\\
  S  \hspace{2.4em} &&  \ \ 1 \hspace{1.5em} &&  \hskip-2mm  \quad \bbR\times (\mathbb{C}^{\times} )^{p-1}\times S^1&&  \hskip-2mm  2 \quad
  \end{align*}

In this table, if $S= \{ \alpha_1, \alpha_3,\cdots, \alpha_{2p-3}, \alpha_{2p-1}  \}$,  $i\in \{1, 3, \dots, 2p-3, 2p-1 \} \text{ and } \{i_1, \dots, i_k\}\subset \{ 1, 3, \dots, 2p-3, 2p-1\}$; 
if $S= \{ \alpha_1, \alpha_3,\cdots, \alpha_{2p-3}, \alpha_{2p-1}  \}$,  $i\in \{1, 3, \dots, 2p-3, 2p-1 \} \text{ and } \{i_1, \dots, i_k\}\subset \{ 1, 3, \dots, 2p-3, 2p-1\}$.

Then from (\ref{e:coeff}), 
$$
K_S= \sum_{k=0} ^{p-1} (-1)^k {p \choose k} +(-1)^p\cdot 2=
\begin{cases}
1& \text{ if $p$ is even,}\\
-1 & \text{ if $p$  is odd}. 
\end{cases}
$$

\item[(iii)] \underline{Type $D_{n}$, $n=2p+1$}

This case is the same as the case of $n=2p$ and $S=\{ \alpha_{n-1}, \alpha_n\}$.

\item[(iv)] \underline{Type $E_7$}

When $S= \{ \alpha_1, \alpha_3, \alpha_7  \}$:
\medskip
\begin{align*}
   S'\hspace{3.2em} &&\hskip-2mm  \#\{ S' \}   \hspace{1em}  &&  H_{S'} \   \  \ &&  \  \ \hskip-2mm  C(H_{S'}) \\
   \emptyset \hspace{3.3em}   &&  \ \ 1 \hspace{1em} \ \  &&\hskip-2mm  (\bbR^\times)^{7}&&  \hskip-2mm  1 \quad \\
 \{\alpha_i\}   \hspace{2.7em}&&  \ \ 3 \hspace{1.7em}  &&  \hskip-2mm  \quad (\mathbb{R}^{\times})  ^{5} \times\mathbb{C}^{\times} &&  \hskip-2mm  1 \quad \\
  \{\alpha_{i_1}, \alpha_{i_2}\}   \hspace{2em}&&  \ \ 3 \hspace{1.7em}  &&  \hskip-2mm  \quad (\mathbb{R}^{\times})  ^{3} \times (\mathbb{C}^{\times} )^2 &&  \hskip-2mm  1 \quad \\
  S  \hspace{3.2em} &&  \ \ 1 \hspace{1.7em} &&  \hskip-2mm  \quad (\bbR^\times)^{n-3}\times \mathbb{C}^{\times} \times S^1&&  \hskip-2mm  2 \quad
  \end{align*}

In this table, $i\in \{ 1, 3, 7 \}$, and $\{i_1, i_2\}\subset \{1, 3, 7\}$. Then from (\ref{e:coeff}), 
$$
K_S 
=1-3+3-2=-1.$$

Since $K_S\neq 0$ for all $J(c_S(Sh))$ with $S\in R_D$, we conclude that $\prod _{R_D}(\tu G)\subseteq $Lift$(\bbC)$.

\end{itemize}

\end{proof}

\end{theorem}

We summarize this paper by the following proposition.

\begin{prop}
Let $\tu G$ be the nonlinear double cover of the simply laced split real group $G$ with $G=SL(n,\bbR), Spin(n,n)$ or the split real group of type $E$. Let $\cal{SH}_{\rho/2} =\{Sh_i\}$ be
the set of irreducible quotients of pseudospherical principal series representations of $\tu G$ with infinitesimal character $\rho/2$.  The roots $\alpha_j$ are labeled according to 
\cite[Table 3]{Ts}. Then we have 
\begin{itemize}
\item[(a)] When $G=SL(2p, \bbR)$, 
\begin{center} \emph{Lift}$_G^{\tu G} (\bbC) = Sh_1+Sh_2+(-1)^p (\pi_1+\pi_2)$, \end{center}
where $\pi_i = c_{\{ \alpha_1,\alpha_3, \dots, \alpha_{2p-1}\}} (Sh_i)$.

\item[(b)] When $G=SL(2p+1, \bbR)$, 
\begin{center}\emph{Lift}$_G^{\tu G} (\bbC) = Sh.$ \end{center}

\item[(c)] When $G=Spin(2p, 2p)$, 
$$
\text{\emph{Lift}}_G^{\tu G} (\bbC) = \sum _{i=1} ^4 \left [ Sh_i+\pi_i + (-1)^p (\delta_i+\sigma_i) \right ],
$$
where $ \pi_i =c_{\{  \alpha_{2p-1}, \alpha_{2p} \}} (Sh_i),  \  \delta_i = c_{\{ \alpha_1,\alpha_3, \dots, \alpha_{2p-3},  \alpha_{2p-1}\}} (Sh_i), \ \delta_i = c_{\{ \alpha_1,\alpha_3, \dots,  \alpha_{2p-3}, \alpha_{2p}\}} (Sh_i)$.

\item[(d)] When $G=Spin(2p+1, 2p+1)$, 
$$
\text{\emph{Lift}}_G^{\tu G} (\bbC) = Sh_1+Sh_2+\pi_1+\pi_2,
$$
where $ \pi_i =c_{\{  \alpha_{2p-1}, \alpha_{2p} \}} (Sh_i)$.

\item[(e)] When $G=E_6$\emph{(split)}, 
\begin{center} \emph{Lift}$_G^{\tu G} (\bbC) = Sh.$ \end{center}

\item[(f)] When $G=E_7$\emph{(split)}, 
$$
\text{\emph{Lift}}_G^{\tu G} (\bbC) = Sh_1+Sh_2-\pi_1-\pi_2,
$$
where $ \pi_i =c_{\{  \alpha_1, \alpha_3, \alpha_7 \}} (Sh_i)$.

\item[(g)] When $G=E_8$\emph{(split)}, 
\begin{center} \emph{Lift}$_G^{\tu G} (\bbC) = Sh.$ \end{center}

\end{itemize}

\end{prop}

\end{document}